\documentclass[11pt,a4paper]{article}
\usepackage{amsmath,amssymb,latexsym,amsthm}
\usepackage[english]{babel}
\usepackage[latin2]{inputenc}

\newtheorem{thm}{Theorem}
 \newtheorem{Coro}[thm]{Corollary}
 \newtheorem{Lemma}[thm]{Lemma}
 \newtheorem{Prop}[thm]{Proposition}
 \theoremstyle{definition}
 
 \theoremstyle{remark}
 
 \newcommand{\Pro}{\noindent{\em Proof. }}

\begin{document}

\centerline{\Large \bf Three-dimensional topological loops with solvable}
\centerline{\Large \bf multiplication groups}

\begin{abstract}
We prove that each $3$-dimensional connected topological loop $L$ having a solvable Lie group of dimension $\le 5$ as the multiplication group of $L$
is centrally nilpotent of class $2$. Moreover, we classify the solvable non-nilpotent Lie groups $G$  which are multiplication groups for $3$-dimensional simply connected topological loops $L$ and $\hbox{dim} \ G \le 5$. These groups are direct products of proper connected Lie groups and have dimension $5$. We find also the inner mapping groups of $L$.
\end{abstract}



\bigskip
\centerline{\bf 1. Introduction}

\medskip
\noindent
The multiplication group $Mult(L)$ of a loop $L$
introduced in \cite{albert1}, \cite{bruck} connects the loop with the group theory since for any normal subloop of $L$ there is a normal subgroup of $Mult(L)$ and conversely to every normal subgroup of $Mult(L)$ corresponds a normal subloop of $L$ (cf. Lemma \ref{brucklemma}).
Necessary and sufficient conditions for a group $K$ to be the multiplication group $Mult(L)$ of a loop $L$ are established in \cite{kepka}. In this criterion there are two special transversals $A$ and $B$ with respect to a subgroup $S$ (see Lemma \ref{kepka}) which results in being the stabilizer of the identity of $L$ in $Mult(L)$ and it is called the inner mapping group $Inn(L)$ of $L$. For finite loops the importance of $Mult(L)$ and $Inn(L)$ as well as the transversals $A$ and $B$ is documented in (\cite{niem3}, \cite{niem4} - \cite{niem2}, \cite{vesanen}).

In general the multiplication group $Mult(L)$ for a topological loop $L$ has infinite dimension.
 If $L$ has a Lie group as its multiplication group, then the structure of $L$ as well as that of $Mult(L)$ is strongly restricted. Hence it is justified to investigate Lie groups which are multiplication groups of $L$ (\cite{figula0} - \cite{figula2}, \cite{loops}). In this case
the criterion in \cite{kepka} can be effectively used and the topological loop $L$ is realized as a sharply transitive section in a subgroup $G$ of $Mult(L)$. This subgroup $G$ is the group topologically generated by the left translations of $L$.

If the group $Mult(L)$ of a $2$-dimensional topological loop $L$ is a Lie group, then it is an elementary filiform Lie group ${\mathcal F}_n$ with $n \ge 4$  (\cite{figula0}).
Classifying all at most $5$-dimensional solvable non-nilpotent Lie groups $K$ which are multiplication groups $Mult(L)$ of $3$-dimensional connected simply connected topological loops $L$ we see that for the structure of $Mult(L)$ one has more freedom. Moreover, knowing $Mult(L)$
one can describe the structure of $L$ and determine the inner mapping group of $L$.

In Section 3 we give the precise structure of the $3$-dimensional simply connected topological loops $L$ such that the multiplication group $Mult(L)$ of $L$ is a solvable Lie group and $L$ has a $1$-dimensional connected normal subloop (see Theorem \ref{solvableonedimensional}).
In this paper 
we prove that for each $3$-dimensional simply connected topological loop $L$ having a solvable Lie group
of dimension $\le 5$  as the multiplication group $Mult(L)$ of $L$ the group $Mult(L)$ is a semidirect product of a group $Q \cong \mathbb R^2$ with 
the group $M= Z \times Inn(L) \cong \mathbb R^n$, $n \in \{ 2,3 \}$, where $\mathbb R=Z$ is a central subgroup of $Mult(L)$. So we show that none of the $4$-dimensional solvable Lie groups as well as none of the $5$-dimensional solvable non-nilpotent indecomposable Lie groups
are multiplication groups of $3$-dimensional topological loops $L$ (see Sections 4 and 5). But there are many loops $L$ having a $4$-dimensional solvable Lie group as the group generated by their left translations (Theorem \ref{lefttranslation}).

To classify the $5$-dimensional solvable decomposable Lie groups $Mult(L)$ of $L$ 
 we have to find special left transversals to a $2$-dimensional subgroup $S$ of $Mult(L)$ such that the core of $S$ in $Mult(L)$ is trivial, $S$ is included in a normal subgroup $M \cong \mathbb R^3$ of $Mult(L)$ with $Mult(L)/M \cong \mathbb R^2$ and the normalizer of $S$ in $Mult(L)$
is the direct product of $S$ and the centre of $Mult(L)$.
The final result of our efforts is the following: If $Mult(L)$ has $1$-dimensional centre, then it is either the group
${\mathcal F}_3 \times {\mathcal L}_2$ or the group $\mathbb R \times {\mathcal L}_2 \times {\mathcal L}_2$, or the direct product
$\mathbb R \times \Sigma$, where $\Sigma$ is a $4$-dimensional indecomposable solvable Lie group having $2$-dimensional commutator subgroup and at most one $1$-dimensional normal subgroup. If $Mult(L)$ has $2$-dimensional centre, then $Mult(L)$ is either the group
${\mathcal F}_{4} \times \mathbb R$ or the direct product of $\mathbb R^2$ and a
$3$-dimensional Lie group having $2$-dimensional commutator subgroup (see Theorem \ref{multinn}).

We want to mention that a Lie group  need not to be the multiplication group of a topological loop if its universal covering has this property. We illustrate this for the direct product $\Omega $ of $\mathbb R^2$ and the group of orientation preserving motions of the euclidean plane and the universal covering of $\Omega $ (Theorem \ref{multinn} case 6) and Proposition \ref{sincos3dim}).

As our result did not give any example of a $3$-dimensional topological loop $L$ having an indecomposable solvable Lie group as the multiplication group of $L$,  further investigations should be focused on this type of groups.

\bigskip
\noindent
\centerline{\bf 2. Preliminaries}

\medskip
\noindent
A binary system $(L, \cdot )$ is called a loop if there exists an element
$e \in L$ such that $x=e \cdot x=x \cdot e$ holds for all $x \in L$ and the
equations $a \cdot y=b$ and $x \cdot a=b$ have precisely one solution, which we denote by $y=a \backslash b$ and $x=b/a$. A loop $L$ is proper if it is not a group.

The left and right translations $\lambda _a: y \mapsto a \cdot y:L \times L \to L$ and $\rho _a: y \mapsto y \cdot a: L \times L \to L$, $a \in L$, are bijections of $L$. The permutation group $Mult(L)$ generated by all left and right translations of the loop $L$ is called the multiplication group of $L$ and the stabilizer of $e \in L$ in the group $Mult(L)$ is called the inner mapping group $Inn(L)$ of $L$.

Let $K$ be a group, let $S \le K$, and let $A$ and $B$ be two left transversals to $S$ in $K$. We say that $A$ and $B$ are $S$-connected if $a^{-1} b^{-1} a b \in S$ for every $a \in A$ and $b \in B$. The core $Co_K(S)$ of $S$ in $K$ is the largest normal subgroup of $K$ contained in $S$. If $L$ is a loop, then
$\Lambda (L)=\{ \lambda _a; \ a \in L \}$ and $R (L)=\{ \rho _a; \ a \in L \}$ are $Inn(L)$-connected transversals in the group $Mult(L)$, and the core of $Inn(L)$ in $Mult(L)$ is trivial. We often use the following (see  \cite{kepka}, Theorem 4.1 and Proposition 2.7).

\begin{Lemma} \label{kepka} A group $K$ is isomorphic to the multiplication group of a loop if and only if there exists a subgroup $S$ with
$Co_K(S)=1$ and $S$-connected transversals $A$ and $B$ satisfying $K=\langle A, B \rangle $.
\end{Lemma}

\begin{Lemma} \label{niemenmaa} Let $L$ be a loop with multiplication group $Mult(L)$ and inner mapping group $Inn(L)$. Then
the normalizer $N_{Mult(L)}(Inn(L))$ is the direct product $Inn(L) \times Z(Mult(L))$, where $Z(Mult(L))$ is the centre of the group $Mult(L)$.
\end{Lemma}

\medskip
\noindent
The kernel of a homomorphism $\alpha :(L, \cdot ) \to (L', \ast )$ of a loop $L$ into a loop $L'$ is a normal subloop $N$ of $L$.
The centre $Z(L)$ of a loop $L$ consists of all elements $z$ which satisfy the equations $z x \cdot y=z \cdot x y, \ x \cdot y z=x y \cdot z, \ x z \cdot y=x \cdot z y, \ z x =x z$ for  all $x,y \in L$.
If we put $Z_0=e$, $Z_1=Z(L)$ and $Z_i/Z_{i-1}=Z(L/Z_{i-1})$, then we obtain a series of normal subloops of $L$. If $Z_{n-1}$ is a proper subloop of $L$ but $Z_n=L$, then $L$ is centrally nilpotent of class $n$. The next assertion was proved by Albert in  \cite{albert1}, Theorems 3, 4 and 5 and by Bruck in \cite{bruck}, IV.1, Lemma 1.3.

\begin{Lemma} \label{brucklemma} Let $L$ be a loop with multiplication group $Mult(L)$ and identity element $e$.
\newline
\noindent
(i) Let $\alpha $ be a homomorphism of the loop $L$ onto the loop $\alpha (L)$ with kernel $N$. Then $\alpha $ induces a homomorphism of the group $Mult(L)$ onto the group $Mult(\alpha (L))$.

Let $M(N)$ be the set $\{ m \in Mult(L); \ x N=m(x) N \ \hbox{for  all} \ x \in L \}$.
Then $M(N)$ is a normal subgroup of $Mult(L)$ containing the multiplication group $Mult(N)$ of the loop $N$ and the multiplication group of the factor loop $L/N$ is isomorphic to $Mult(L)/M(N)$.
\newline
\noindent
(ii) For every normal subgroup $\mathcal{N}$ of $Mult(L)$ the orbit $\mathcal{N}(e)$ is a normal subloop of $L$. Moreover,
$\mathcal{N} \le  M(\mathcal{N}(e))$.
\end{Lemma}

\medskip
\noindent
A loop $L$ is called topological  if $L$ is a topological space and the binary operations
$(x,y) \mapsto x \cdot y, \ (x,y) \mapsto x
\backslash y, (x,y) \mapsto y/x :L \times L \to L$  are continuous.
Let $G$ be a connected Lie group, let $H$ be a subgroup of $G$.  A continuous section $\sigma :G/H \to G$ is called sharply transitive, if the set $\sigma (G/H)$ operates sharply transitively on $G/H$, which means that for any $x H$ and $y H$ there exists precisely one
$z \in \sigma (G/H)$ with $z x H= y H$.
Every connected topological loop $L$ having a Lie group $G$ as the group topologically generated by the left translations of $L$ is obtained on a homogeneous space $G/H$, where $H$ is a closed subgroup of $G$ with $Co_G(H)=1$ and $\sigma :G/H \to G$ is  a continuous sharply transitive section such that $\sigma (H)=1 \in G$ and the subset $\sigma (G/H)$ generates $G$.  The multiplication  of  $L$ on the manifold  $G/H$  is  defined by  $x H \ast y H=\sigma (x H) y H$ and  the group $G$ is  the group topologically generated by the left translations of $L$. Moreover, the subgroup $H$ is the stabilizer of the identity element $e \in L$ in the group $G$. The following assertion is proved in \cite{hofmann}, IX.1.

\begin{Lemma} \label{simplyconnecteduj} For any connected topological loop there is a universal covering loop. This loop  is simply connected. \end{Lemma}

\noindent
The elementary filiform Lie group ${\mathcal F}_{n}$ is the
 simply connected Lie group of dimension $n \ge 3$ such that its Lie algebra has a basis $\{ e_1, \cdots , e_{n} \}$  with  $[e_1, e_i]= e_{i+1}$ for
$2 \le i \le n-1$.  A $2$-dimensional simply connected loop $L_{\mathcal F}$ is called
an elementary filiform loop if its multiplication group is an elementary filiform group  ${\mathcal F}_{n}$, $n \ge 4$ (\cite{figula}).

Homogeneous spaces of solvable Lie groups are called solvmanifolds.

\bigskip
\noindent
\centerline{\bf 3. Three-dimensional topological loops with}
\centerline{\bf one-dimensional connected normal subloop}

\medskip
\noindent
Let $L$ be a topological loop on a connected $3$-dimensional manifold such that the group $Mult(L)$ topologically generated by all left and right  translations of $L$ is a Lie group. The loop $L$ is a $3$-dimensional homogeneous space with respect to the transformation group $Mult(L)$ acting transitively and effectively on $L$. According to Theorem B and Theorem 1 in \cite{mostow} the simply connected spaces $S^2 \times \mathbb R$ and $S^3$ are not solvmanifolds. Hence from \cite{gorbatsevich} we get the following.

\begin{Lemma} \label{simplyconnected} Let $L$ be a $3$-dimensional proper connected topological loop such that its multiplication group $Mult(L)$ is a solvable Lie group. If $L$ is simply connected, then it is homeomorphic to $\mathbb R^3$.
\end{Lemma}

\medskip
\noindent
Assume that the multiplication group $Mult(L)$ of a topological loop $L$ is solvable. Let $K$ be a minimal non-trivial connected normal subgroup of $Mult(L)$. Then one has $\hbox{dim} \ K \in \{1, 2\}$. By Lemma \ref{brucklemma} the orbit $K(e)$ is a connected normal subloop of $L$.
Since the core $Co_{Mult(L)}(Inn(L))$ is trivial $K(e) \neq \{ e \}$. Hence the dimension of $K(e)$ is $1$ or $2$. Now we deal with the case that $\hbox{dim} \ K(e)=1$.

\begin{thm} \label{solvableonedimensional} Let $L$ be a $3$-dimensional proper connected simply connected topological loop such that its multiplication group $Mult(L)$ is a solvable Lie group. If $L$ has a $1$-dimensional connected normal subloop $N$, then $N$ is isomorphic to the group $\mathbb R$ and we have the following possibilities:
\newline
\noindent
(a) The factor loop $L/N$ is isomorphic to $\mathbb R^2$. Then $N$ is contained in the centre of $L$ and
the group $Mult(L)$ is a semidirect product of a group $Q \cong \mathbb R^2$ with 
the abelian group $M= Z \times Inn(L) \cong \mathbb R^m$,
$m \ge 2$, where $\mathbb R=Z \cong N$ is a central subgroup of $Mult(L)$.
\newline
\noindent
(b) The loop $L/N$ is isomorphic either to the non-abelian $2$-dimensional Lie group
${\mathcal L}_2$ or to a $2$-dimensional elementary filiform loop $L_{\mathcal F}$. Then the group $Mult(L)$ has a normal subgroup $S$ containing
$Mult(N) \cong \mathbb R$ such that the factor group $Mult(L)/S$ is isomorphic to the direct product ${\mathcal L}_2 \times {\mathcal L}_2$ if
$L/N \cong {\mathcal L}_2$ or to an elementary filiform Lie group ${\mathcal F}_{n}$, $n \ge 4$, if $L/N \cong L_{\mathcal F}$. Moreover, $Mult(L)$ has dimension at least $5$.
\end{thm}
\Pro By Lemma \ref{simplyconnected} the loop $L$ is homeomorphic to $\mathbb R^3$. The connected normal subloop $N$ of $L$ is isomorphic to
$\mathbb R$ because 
 the multiplication group of $N$ a Lie subgroup of $Mult(L)$ (Theorem 18.18 in \cite{loops}). The factor loop $L/N$ is a $2$-dimensional connected loop such that the multiplication group $Mult(L/N)$ is a factor group of $Mult(L)$ (Lemma \ref{brucklemma}). The manifold $L$ is a fibering of $\mathbb R^3$ over $L/N$ with fibers homeomorphic to $\mathbb R$. Hence $L/N$ is homeomorphic to $\mathbb R^2$
and therefore it is either a $2$-dimensional connected Lie group or an elementary filiform loop $L_{\mathcal F}$ (Theorem 1 in \cite{figula0}).

If the factor loop $L/N$ is the Lie group $\mathbb R^2$, then by Lemma \ref{brucklemma} there exists a normal subgroup $M$ of $Mult(L)$ such that $Mult(L)/M$ is isomorphic to the multiplication group of the loop $L/N$ and hence to the group $\mathbb R^2$. Therefore the group $M$ is connected and $Mult(L)/M$ operates sharply transitively on the orbits of $N$ in $L$. The group $M$ contains the multiplication group $Mult(N) \cong \mathbb R$ of $N$ and leaves every orbit of $N$ in the manifold $L$ invariant. Every orbit of $N$ is homeomorphic to
$\mathbb R$. Hence the group $M$ induces on the orbit $N(e)$ either the sharply transitive group $\mathbb R$ or the group $\Omega $ isomorphic to the Lie group ${\mathcal L}_2$ (\cite{salzmann}, Lemma 1.10).

Assume first that the group induced by $M$ on $N(e)$ is $\Omega \cong {\mathcal L}_2$. Then $M$ induces a group isomorphic to $\Omega $ on every orbit $N(x)$, $x \in L$. Since all  $1$-dimensional connected subgroups of $\Omega $ different from the commutator subgroup are conjugate, the stabilizer $\Omega_e$ of $e \in L$ in $\Omega $ would fix on every orbit $N(x)$ precisely one point. The set of fixed points of $\Omega_e$ in $L$ coincides with that of fixed points of the stabilizer $Inn(L)$ of $e \in L$ in $Mult(L)$. This latter is the centre $Z$ of $L$ (see \cite{bruck}, IV.1). Hence the centre $Z$ of $L$ would be at least $2$-dimensional and we would have $L= N \cdot Z$. But then $L$ would be an abelian group which is a contradiction.

Therefore the group $M$ induces on every orbit $N(x)$, $x \in L$, the sharply transitive group $\mathbb R$. The stabilizer
$M_1$ of $e \in L$ in $M$ fixes every point of the orbit $N(e)=M(e)$. Hence $M_1$ is a normal subgroup of $M$. Since the factor group $M/M_1$ is isomorphic to  $\mathbb R$ the commutator subgroup $M'$ of $M$ is contained in $M_1$ and $M'$ is normal in $Mult(L)$. If $M'$ were different from
$\{ 1 \}$, then  $Mult(L)$ would contain the normal subgroup $M'$ which has fixed points. This is a contradiction because the transitive group $Mult(L)$ acts effectively on $L$. Hence $M$ is abelian. If $M$ would contain a compact connected subgroup $K \neq \{ 1 \}$, then $K$ would be isomorphic to the group $SO_2(\mathbb R)$  and it would be a normal subgroup of $Mult(L)$ which has a fixed point in $L$. This contradiction yields that $M$ is isomorphic to $\mathbb R^n$. Since $L$ is a proper loop of dimension $3$ one has $\hbox{dim} \ Mult(L) \ge 4$ and hence $n \ge 2$.   
As the inner mapping group $Inn(L)$ has codimension $3$ it is the group $M_1$.
Since $M_1$ fixes every element of the loop $N(e)$ the normal subloop $N$ is a central subgroup of $L$. The group consisting of the translations by elements of $N$ is isomorphic to $N$ and it is a central subgroup $Z$ of $Mult(L)$. Then we have $M=Z \times Inn(L)$ and the assertion (a) is proved.

If the factor loop $L/N$ is isomorphic to the Lie group ${\mathcal L}_2$, respectively to an elementary filiform loop $L_{\mathcal F}$, then the multiplication group $Mult(L/N)$ is isomorphic to the direct product
${\mathcal L}_2 \times {\mathcal L}_2$, respectively to an  elementary filiform Lie group ${\mathcal F}_{n}$, $n \ge 4$. Moreover, there exists a normal subgroup $S$ of $Mult(L)$ containing the group $Mult(N) \cong \mathbb R$
(see Lemma \ref{brucklemma}) such that $Mult(L)/S$ is isomorphic to the group $Mult(L/N)$ and the assertion (b) follows.  \qed

\bigskip
\noindent
\centerline{\bf 4. Three-dimensional topological loops with four-dimensional}
\centerline{\bf  solvable Lie group as multiplication group do not exist}

\medskip
\noindent
The following Lemma follows from Theorem 18.18 in \cite{loops}, Theorem 1 in \cite{figula0} and Theorem \ref{solvableonedimensional} (a).

\begin{Lemma} \label{fourdim}
If there exists proper connected topological loop $L$ having a $4$-dimensional solvable non-nilpotent Lie group as its multiplication group $Mult(L)$, then $L$ has dimension $3$. Moreover, if $L$ is simply connected and has a $1$-dimensional normal subloop, then $Mult(L)$ is a semidirect product of $\mathbb R^2$ with a normal subgroup $M \cong \mathbb R^2$ containing a $1$-dimensional central subgroup of $Mult(L)$. 
\end{Lemma}

\medskip
\noindent
The $4$-dimensional indecomposable Lie algebras are listed in \cite{mubarakjzanov1}, § 5.
Among these solvable Lie algebras there are four with $1$-dimensional centre: the filiform Lie algebra ${\it g}_{4,1}$ and the non-nilpotent Lie algebras ${\it g}_{4,3}$, ${\it g}_{4,8}$ with $h=-1$, ${\it g}_{4,9}$ with $p=0$.  Proposition 4.3 in \cite{figula} shows that the filiform Lie group ${\mathcal F}_{4}$ is not the multiplication group of $3$-dimensional connected topological loops. Since the commutator Lie algebra of
${\it g}_{4,8}$ and ${\it g}_{4,9}$ has dimension $3$ there is no connected topological loop $L$ having these Lie algebras as the Lie algebra
of $Mult(L)$ (see Lemma \ref{fourdim} and Theorem \ref{solvableonedimensional} (a)).

The commutator Lie algebra of ${\it g}_{4,3}$ has dimension $2$. Hence for the corresponding simply connected Lie group $G$ it seems to be more natural that $G$ can be the multiplication group $Mult(L)$ of connected topological loops. Although, as we will show, there are four classes of  $3$-dimensional simply connected topological loops $L$ having $G$ as the group generated by their left translations (Theorem \ref{lefttranslation}), for any of these loops the multiplication group $Mult(L)$ has dimension greater than $4$ (Corollary \ref{corro}). For the classification of these loops $L$ we often use the following lemmata, the first of which is proved in \cite{figula} Lemma 4.2, and the second in \cite{figula2} Lemma 3.1.

\begin{Lemma} \label{bijective}
Let $f: (x,y,z) \mapsto f(x,y,z): \mathbb R^3 \to \mathbb R$ be a continuous function.
The function  $g: z \mapsto z + u f(x_0,y_0,z): \mathbb R \to \mathbb R$ is bijective for every $x_0,y_0,u \in \mathbb R$
if and only if $f$ does not depend on the variable $z$.
\end{Lemma}

\begin{Lemma} \label{functional} Let $f: \mathbb R \to \mathbb R$ be a continuous function such that for all $z_1, z_2 \in \mathbb R$ one has 
$f(z_2)+ e^{-z_2} f(z_1) = f(z_1+z_2)$. 
Then we get $f(z)=c(1-e^{-z})$, where $c$ is a real constant.
\end{Lemma}

\begin{thm} \label{lefttranslation} Let $G$ be the four-dimensional connected simply connected solvable Lie group the multiplication of which is represented on $\mathbb R^4$ by
\begin{equation} \label{multgroup1} g(x_1,x_2,x_3,x_4) g(y_1,y_2,y_3,y_4)= g(x_1+ y_1 e^{x_4}, x_2+y_2+x_4 y_3, x_3+y_3, x_4+y_4). \nonumber \end{equation}
Let $H$ be a non-normal subgroup of $G$ isomorphic to $\mathbb R$. Using suitable automorphisms of $G$ we may choose $H$ as one of the following subgroups:
\[ H_1= \{ g(0,0,0,x_4); x_4 \in \mathbb R \}, \ H_2= \{ g(0,0,x_3,0); x_3 \in \mathbb R \}, \]
\[ H_3=\{ g(x_3,0,x_3,0); x_3 \in \mathbb R \}, \ H_4=\{ g(x_1,x_1,0,0); x_1 \in \mathbb R \}. \]
\noindent
a) Every continuous sharply transitive section $\sigma : G/H_1 \to G$ with the properties that $\sigma (G/H_1)$ generates $G$ and $\sigma (H_1)=1$ is determined by the map $\sigma _f : g(x,y,z,0) H_1 \mapsto g(x,y,z,f(z))$, where $f: \mathbb R \to \mathbb R$ is a continuous non-linear function with
$f(0)=0$. The multiplication of the loop $L_f$ given by $\sigma _f$ can be written as
\begin{equation} \label{fourdimelso} (x_1,y_1,z_1) \ast (x_2,y_2,z_2)= (x_1+ x_2 e^{f(z_1)},y_1+y_2+z_2 f(z_1),z_1+z_2). \end{equation}

\noindent
b) Each continuous sharply transitive section $\sigma : G/H_2 \to G$ such that $\sigma (G/H_2)$ generates $G$ and
$\sigma (H_2)=1$ has the form
\[ \sigma _h: g(x,y,0,z) H_2 \mapsto g(x,y+ h(x,z) z, h(x,z),z), \]
where $h: \mathbb R^2 \to \mathbb R$ is a continuous function with $h(0,0)=0$ such that $h$ does not fulfil the identities $h(x,0)=0$ and $h(0,z)=l z$, $l \in \mathbb R$, simultaneously.
The multiplication of the loop $L_h$ corresponding to $\sigma _h$ is determined by
\begin{equation} \label{multiplicationuj2} (x_1,y_1,z_1) \ast (x_2,y_2,z_2)=  (x_1+ x_2 e^{z_1}, y_1+y_2- z_2 h(x_1,z_1), z_1+z_2). \end{equation}

\noindent
c) Every continuous sharply transitive section $\sigma : G/H_3 \to G$ such that $\sigma (G/H_3)$ generates $G$ and $\sigma (H_3)=1$ is given by the map
\[ \sigma _f: g(x,y,0,z) H_3 \mapsto g(x+ e^{z} f(x,y,z), y+ z f(x,y,z), f(x,y,z),z)  \]
with a continuous function
 $f: \mathbb R^3 \to \mathbb R$ such that $f(0,0,0)=0$, $f$ does not satisfy either the identities
\begin{equation} \label{elsofeltetel} f(x,y,0)=-x, \  f(0,0,z)=C(1-e^{-z}), \ C \in \mathbb R,  \end{equation}
or the identities
\begin{equation} \label{masodikfeltetel} f(x,y,0)=0, \ f(0,0,z)= \lambda z, \ \lambda \in \mathbb R,  \end{equation}
simultaneously and
for all triples $(x_1,y_1,z_1)$ and $(x_2,y_2,z_2) \in \mathbb R^3$ the equations
\begin{equation} \label{a43h3masodikequuj} y= y_2- y_1+ z_1 f(x,y,z_2-z_1),  \end{equation}
\begin{equation} \label{a43h3harmadikequuj} x= x_2 -e^{z_2-z_1}x_1 + e^{z_2} (1-e^{-z_1}) f(x,y,z_2-z_1) \end{equation}
have a unique solution $(x,y) \in \mathbb R^2$.
The loop $L_f$ corresponding to $\sigma _f$ is defined by the multiplication
\begin{equation} \label{multiplicationuj4} (x_1,y_1,z_1) \ast (x_2,y_2,z_2)=  \nonumber \end{equation}
\begin{equation} (x_1+ e^{z_1}(x_2+ f(x_1,y_1,z_1)(1-e^{z_2})), y_1+y_2- z_2 f(x_1,y_1,z_1), z_1+z_2). \end{equation}

\noindent
d) Any continuous sharply transitive section $\sigma : G/H_4 \to G$ such that $\sigma (G/H_4)$ generates $G$ and
$\sigma (H_4)=1$ is determined by the map
\[ \sigma _k: g(x,0,y,z) H_4 \mapsto g(x+e^{z} k(x,y,z), k(x,y,z), y,z), \]
 where $k: \mathbb R^3 \to \mathbb R$ is a continuous function with $k(0,0,0)=0$ such that $k$ does not fulfil the identities given by (\ref{elsofeltetel}) in case c) simultaneously and such that
for all triples $(x_1,y_1,z_1)$ and $(x_2,y_2,z_2) \in \mathbb R^3$ the equation
\begin{equation} \label{equd} x + e^{z_2} k(x,y_2-y_1,z_2-z_1)[e^{-z_1}-1] = x_2- x_1 e^{z_2-z_1}+ e^{z_2}(z_2-z_1) y_1   \end{equation}
has a unique solution $x \in \mathbb R$.
The multiplication of the loop $L_k$ corresponding to $\sigma _k$ can be written as
\begin{equation} (x_1,y_1,z_1) \ast (x_2,y_2,z_2)= \nonumber \end{equation}
\begin{equation} \label{multiplication1harmadik} (x_1+ e^{z_1}[x_2 + k(x_1,y_1,z_1)-e^{z_2}(z_1 y_2+k(x_1,y_1,z_1))], y_1+y_2,z_1+z_2). \end{equation}   \end{thm}
\Pro
The linear representation of the group $G$ is given in \cite{ghanam}, Case 4.3. Let $L$ be a $3$-dimensional connected simply connected topological loop having $G$ as the group topologically generated by its left translations. Then the stabilizer $H$ of $e \in L$ in $G$ is a $1$-dimensional non-normal subgroup of $G$. As the Lie algebra ${\bf g}$ of $G$
 has a basis $\{ e_1, e_2, e_3, e_4 \}$ with $[e_1, e_4]=e_1$,
$[e_3, e_4]=e_2$,
the subgroup $\exp t e_2$, $t \in \mathbb R$, is the centre of $G$, the subgroup $\exp (t e_2 + s e_1)$,
$t, s \in \mathbb R$, is the commutator subgroup of $G$.  Hence the automorphism group
of ${\bf g}$ consists of the following linear mappings $\varphi(e_1)= a e_1$, $\varphi(e_2)= b e_2$, $\varphi(e_3)= k e_2+ b e_3$,
$\varphi(e_4)= l_1 e_1+ l_2 e_2+ l_3 e_3+ e_4$,
with $a b \neq 0$, $k, l_1, l_2, l_3 \in \mathbb R$.
Since $\mathbb R e_1$ and $\mathbb R e_2$ are ideals of ${\bf g}$ the subalgebra ${\bf h}$ of $H$ does not contain $e_1$, $e_2$. Hence $H$ is a subgroup
$\exp t(\alpha e_1 + \beta e_2 + \gamma e_3 + \delta e_4)$ with $t \in \mathbb R$ such that
$\gamma ^2 + \delta ^2 =1$ or $\alpha \beta \neq 0$. Then a suitable automorphism of $G$ corresponding to an automorphism $\varphi $ of
${\bf g}$ maps $H$ onto one of the following subgroups
\[ H_1= \exp{ t e_4}, \ \ H_2= \exp{ t e_3}, \ \ H_3=\exp{ t(e_3+e_1)}, \ \ H_4=\exp{ t(e_1+e_2)}. \]
Every connected topological proper loop $L$ having $G$ as the group topologically generated by its left translations and $H$ as the stabilizer of $e \in L$ in $G$ is determined by a continuous sharply transitive section
$\sigma : G/H \to G$ with the properties that $\sigma(H)=1 \in G$ and $\sigma(G/H)$ generates $G$.

\noindent
First we assume that $H=H_1=\{ g(0,0,0,k); \ k \in \mathbb R \}$. Since all elements of $G$ have a unique decomposition as
$g(x, y, z, 0)  g(0, 0, 0, k)$, any continuous function
$f: \mathbb R^3 \to \mathbb R; (x,y, z) \mapsto f(x,y,z)$ determines a continuous section $\sigma : G/H_1 \to G$ given by
\[ \sigma: g(x,y,z,0) H_1 \mapsto  g(x, y, z, 0) g(0, 0, 0, f(x,y,z)) = g(x,y,z,f(x,y,z)). \]
The section $\sigma $ is sharply transitive if and only if for any triple $(x_1,y_1,z_1)$, $(x_2,y_2,z_2) \in \mathbb R^3$ there exists precisely one triple  $(x,y,z) \in \mathbb R^3$ such that
\[ g(x,y,z,f(x,y,z)) g(x_1, y_1, z_1, 0)= g(x_2, y_2, z_2, 0) g(0,0,0,t) \]
for a suitable $t \in \mathbb R$. This provides the following equations
$z= z_2-z_1$, $t=f(x,y,z_2-z_1)$,
\begin{equation} \label{equa43elso} y= y_2-y_1 -z_1 f(x,y,z_2-z_1),  \end{equation}
\begin{equation} \label{equa43masodik} x= x_2 -x_1 e^{f(x,y,z_2-z_1)}.  \end{equation}
For $x_1=0$ equation (\ref{equa43masodik}) yields that $x=x_2$ and equation (\ref{equa43elso}) has a unique solution for $y$ if and only if the function $g: y \mapsto y+ z_1 f(x_0,y,z_0): \mathbb R \to \mathbb R$ is bijective for every $x_0=x_2$, $z_0=z_2-z_1$ and $z_1 \in \mathbb R$. This is the case precisely if
the function $f(x, y, z)=f(x,z)$ does not depend on the variable $y$ (Lemma \ref{bijective}).
Using this, equations (\ref{equa43elso}) and (\ref{equa43masodik}) are reduced to
\begin{equation} \label{equa43harmadik} y= y_2-y_1 -z_1 f(x,z_2-z_1),  \end{equation}
\begin{equation} \label{equa43negyedik} x= x_2 -x_1 e^{f(x,z_2-z_1)}.  \end{equation}
Applying Lemma \ref{bijective} for the function $e^{f(x,z)}: \mathbb R^2 \to \mathbb R$ we obtain that
equation (\ref{equa43negyedik}) has a unique solution for $x$ precisely if
the function $f(x,z)=f(z)$ does not depend on $x$.
Since in this case equation (\ref{equa43harmadik}) has a unique solution $y=y_2-y_1-z_1 f(z_2-z_1)$ each continuous function
$f: \mathbb R \to \mathbb R$ with $f(0)=0$ defines a loop $L_f$. This loop is proper if $\sigma (G/H_1)$ generates $G$.
The set $\sigma (G/H_1)= \{g(x,y,z,f(z)); x,y,z \in \mathbb R \}$ contains the commutator subgroup
$G'=\{ g(x, y, 0, 0);  x,y \in \mathbb R \}$
and the set
$F=\{ g(0,0, z, f(z));  z \in \mathbb R \}$.
We have $G' \cap F= \{ 1 \}$. Therefore $\sigma (G/H_1)$ does not generate $G$ if the set $F G'/G'$ is a one-parameter subgroup of $G/G'$.
As \[g(\mathbb R, \mathbb R, z_1, f(z_1)) g(\mathbb R, \mathbb R, z_2, f(z_2))= g(\mathbb R, \mathbb R, z_1+z_2, f(z_1)+f(z_2)) \]
this is the case precisely if $f(z)= l z$, $l \in \mathbb R$.
Hence for every non-linear function $f$ there is a topological proper loop $L_f$.
\newline
\noindent
In the coordinate system $(x,y,z) \mapsto g(x,y,z,0)H_1$ the multiplication of $L_f$ is determined if we apply
$\sigma (g(x_1,y_1,z_1,0)H_1)=g(x_1,y_1,z_1,f(z_1))$ to the left coset
$g(x_2,y_2,z_2,0)H_1$ and find in the image coset the element of $G$ which lies in the set $\{ g(x,y,z,0)H_1; \ x,y,z \in \mathbb R \}$. A direct computation yields multiplication (\ref{fourdimelso}) and assertion a) is proved.

\medskip
\noindent
A similar consideration as in the previous case yields that for $H=H_2=\{ g(0,0,k,0); \ k \in \mathbb R \}$
an arbitrary continuous section  $\sigma _2: G/H_2 \to G$ may be given by $\sigma _2: g(x,y,0,z) H_2 \mapsto $
\begin{equation} \label{equu}  g(x, y, 0, z) g(0, 0, h(x,y,z), 0) = g(x,y+z h(x,y,z), h(x,y,z),z), \end{equation}
for $H=H_3=\{ g(t,0,t,0); \ t \in \mathbb R \}$ a continuous section $\sigma _3: G/H_3 \to G$ can be given by
\begin{equation} \label{equuu} \sigma _3: g(x,y,0,z) H_3 \mapsto  g(x, y, 0, z) g(f(x,y,z), 0, f(x,y,z), 0) = \nonumber \end{equation}
\begin{equation} g(x+e^{z} f(x,y,z),y+z f(x,y,z), f(x,y,z),z),  \end{equation}
and for $H=H_4=\{ g(t,t,0,0); \ t \in \mathbb R \}$ a continuous section $\sigma _4: G/H_4 \to G$ may be given by
\begin{equation} \label{equuuu}  \sigma _4: g(x,0,y,z) H_4 \mapsto  g(x,0,y,z) g(k(x,y,z), k(x,y,z), 0, 0) = \nonumber \end{equation}
\begin{equation} g(x+e^{z} k(x,y,z), k(x,y,z),y,z), \end{equation}
where $h: \mathbb R^3 \to \mathbb R$, $f: \mathbb R^3 \to \mathbb R$, $k: \mathbb R^3 \to \mathbb R$ are continuous functions. These sections
$\sigma _i$, $i=2,3,4$, have the property $\sigma _i(H)=1 \in G$ precisely if $h(0,0,0)= f(0,0,0)= k(0,0,0)=0$.

\noindent
The set $\sigma _i (G/H_i)$ given by (\ref{equu}), (\ref{equuu}), (\ref{equuuu}) acts sharply transitively on $G/H_i$ if and only if  for $i=2$  the equation
\begin{equation} \label{elsoequuj} g(x,y+z h(x,y,z), h(x,y,z),z) g(x_1, y_1, 0, z_1)= g(x_2, y_2, 0, z_2) g(0,0,t,0),  \end{equation}
for $i=3$ the equation
\begin{equation} \label{a43h3elsoequ} g(x+e^{z} f(x,y,z),y+z f(x,y,z), f(x,y,z),z) g(x_1, y_1, 0, z_1)= \nonumber \end{equation}
\begin{equation} g(x_2, y_2, 0, z_2) g(t,0,t,0), \end{equation}
for $i=4$ the equation
\begin{equation} \label{a43h4equ} g(x+e^{z} k(x,y,z),k(x,y,z),y,z) g(x_1,0,y_1,z_1)= g(x_2,0, y_2,z_2) g(t,t,0,0) \end{equation}
has a unique solution $(x,y,z) \in \mathbb R^3$ with a suitable $t \in \mathbb R$ for any given triple $(x_1, y_1, z_1)$,
$(x_2,y_2,z_2) \in \mathbb R^3$.
Equation (\ref{elsoequuj}) is equivalent to the following $z= z_2-z_1$, $t= h(x,y,z_2-z_1)$, $x= x_2- e^{z_2-z_1} x_1$ and
\begin{equation} 0=y -y_2+ y_1- z_1 h(x_2- e^{z_2-z_1} x_1,y, z_2-z_1). \nonumber \end{equation}
This last equation has a unique solution for $y$ precisely if the function $h(x, y, z)=h(x,z)$ does not depend on the variable $y$ (cf. Lemma \ref{bijective}).
Equation (\ref{a43h3elsoequ}) yields $z= z_2-z_1$, $t= f(x,y,z_2-z_1)$ and that equations
(\ref{a43h3masodikequuj}), (\ref{a43h3harmadikequuj}) in assertion c) have a unique solution $(x,y) \in \mathbb R^2$.
Moreover, equation (\ref{a43h4equ}) gives $z= z_2-z_1$, $y= y_2-y_1$, $t= y_1(z_2-z_1) + k(x,y_2-y_1,z_2-z_1)$ and that equation (\ref{equd}) in assertion d) has a unique solution $x \in \mathbb R$.

\noindent
Now we investigate under which circumstances the set $\sigma _i (G/H_i)$, $i=2,3,4$, generates the group $G$.

\noindent
The set $\sigma _2 (G/H_2)= \{ g(x,y+z h(x,z), h(x,z),z); x,y,z \in \mathbb R \}$ contains the subgroup
$K_2=\{ g(x, y, h(x,0), 0); \ x,y \in \mathbb R \}$
and the subset $F_2=\{ g(0,z h(0,z), h(0,z),z); \ z \in \mathbb R \}$.
The set $\sigma _3 (G/H_3)$ given by (\ref{equuu}) includes the subgroup
$K_3=\{ g(x+f(x,y,0), y, f(x,y,0), 0); \ x,y \in \mathbb R \}$,
and the subset
$F_3=\{ g(e^{z} f(0,0,z),z f(0,0,z), f(0,0,z),z); \ z \in \mathbb R \}$.
The set $\sigma _4 (G/H_4)$ given by (\ref{equuuu}) contains the subgroup
$K_4=\{ g(x+k(x,y,0), k(x,y,0), y, 0); \ x, \\ y \in \mathbb R \}$
and the subset
$F_4=\{ g(e^{z} k(0,0,z), k(0,0,z), 0,z); \ z \in \mathbb R \}$.
As for all these cases we have $K_i \cap F_i=  \{ 1 \}$ the set $\sigma _i(G/H_i)$, $i=2,3,4$, does not generate $G$ if the group $K_i$ has dimension $2$, for all $h \in F_i$ one has $h^{-1} K_i h =K_i$
and $F_i K_i/K_i$ is a one-parameter subgroup of $G/K_i$.

\noindent
First we consider the pair $(K_2, F_2)$.
 The group $K_2$ has dimension $2$ if the subgroup $\{ g(x, 0, h(x,0), 0); x \in \mathbb R \}$ is a one-parameter subgroup. This is the case  precisely if $h(x,0)=b x$, $b \in \mathbb R$.  For $h= g(0,z h(0,z), h(0,z),z) \in F_2$, $z \neq 0$ we get
$h^{-1} g(x, y, b x, 0) h= g(x e^{-z}, y -b z x, b x,0)$ is an element of $K_2$ if and only if $b=0$. Then the group $K_2$ coincides with the commutator subgroup $G'$ of $G$.
The set $(F_2 G')/G'$ is a one-parameter subgroup precisely if $h(0,z)=l z$, $l \in \mathbb R$.
Therefore any function $h: \mathbb R^2 \to \mathbb R$, which does not satisfy the identities $h(x,0)=0$ and
$h(0,z)= l z$, $l \in \mathbb R$, simultaneously determines a proper topological loop $L_h$. A direct computation yields that the multiplication of $L_h$ corresponding to the section $\sigma _2$ in the coordinate system $(x,y,z) \mapsto g(x,y,0,z) H_2$ is given by (\ref{multiplicationuj2}). This proves assertion b).

\noindent
Now we deal with the pair $(K_3,F_3)$.
The group $K_3$ has dimension $2$ if and only if $f(x,y,0)=c x+d y$, $c,d \in \mathbb R$.
 For $h \in F_3$ with $z \neq 0$ we have
\[ h^{-1} g(x+c x+d y, y, c x+d y, 0) h= g([(c+1)x+dy] e^{-z}, y -z (c x+d y), c x+d y,0). \]
Hence $h^{-1} K_3 h =K_3$ if and only if one has either $c=-1$, $d=0$ or $c=d=0$.

\noindent
In the first case $K_3$ is the normal subgroup  $\widetilde{G}= \{g(0,y,-x,0); \ x,y \in \mathbb R \}$ of $G$, in the second case
$K_3=G'$.
Since
\[ g(e^{z_1} f(0,0,z_1), \mathbb R, \mathbb R, z_1)g(e^{z_2} f(0,0,z_2), \mathbb R, \mathbb R, z_2)= \]
\[ g(e^{z_1+z_2} f(0,0,z_2)+e^{z_1} f(0,0,z_1), \mathbb R, \mathbb R, z_1+z_2) \]
the set $F_3 \widetilde{G}/ \widetilde{G}$ is a one-parameter subgroup of $G/\widetilde{G}$ if and only if for all $z_1$, $z_2 \in \mathbb R$ the identity $f(0,0,z_2)+ e^{-z_2} f(0,0,z_1)= f(0,0,z_1+z_2)$ holds. By Lemma \ref{functional} we obtain
$f(0,0,z)=C(1-e^{-z})$ with $C \in \mathbb R$.
The set $F_3 G'/G'$ is a one-parameter subgroup of $G/ G'$ if and only if one has
 $f(0,0,z)= \lambda z$ for some $\lambda \in \mathbb R$.
The set $\sigma _3(G/H_3)$ does not
generate $G$ if the function $f(x,y,z)$ satisfies either the identities given by (\ref{elsofeltetel}) or the identities given by (\ref{masodikfeltetel}) in assertion c).
A direct computation yields that the multiplication of the loop $L_f$ corresponding to the section $\sigma _3$ in the coordinate system $(x,y,z) \mapsto g(x,y,0,z) H_3$ is given by (\ref{multiplicationuj4}) and the assertion c) is proved.

\noindent
Finally we consider the pair $(K_4,F_4)$.
The group $K_4$ has dimension $2$ if and only if
 $k(x,y,0)=a x+b y$, $a,b \in \mathbb R$.
 For $h \in F_4$, $z \neq 0$ we have
\[ h^{-1} g(x+a x+b y, a x+b y, y, 0) h= g([(a+1)x + by] e^{-z}, -z y + a x+b y, y, 0). \]
Hence we obtain $h^{-1} K_4 h =K_4$ if and only if $a=-1$ and $b=0$. Then the group $K_4$ coincides with the group $\widetilde{G}$ introduced in the previous case. Hence the same consideration as there proves that
the set $\sigma _4(G/H_4)$ does not
generate $G$ if the function $k(x,y,z)$ satisfies the identities given by (\ref{elsofeltetel}).
A direct computation gives that in the coordinate system $(x,y,z) \mapsto $$ g(x,0,y,z)H_4$ the multiplication of the loop $L_k$  is given by (\ref{multiplication1harmadik}) and the assertion d) is proved. \qed

\begin{Coro} \label{corro}
There is no connected topological loop $L$ such that the multiplication group of $L$ is locally isomorphic to the group  $G$ in Theorem \ref{lefttranslation}.  \end{Coro}
\Pro By Lemmata \ref{simplyconnecteduj}, \ref{fourdim}, \ref{simplyconnected} we may assume that $L$ is homeomorphic to $\mathbb R^3$. Every Lie group locally isomorphic to the group $G$ in Theorem \ref{lefttranslation} has a $1$-dimensional centre $Z$. The orbit $Z(e)$ is a $1$-dimensional normal subloop of $L$ isomorphic to $\mathbb R$ (see Lemma \ref{brucklemma}). Hence the multiplication group of $L$ is the simply connected group $G$ (cf. Lemma \ref{fourdim}) and the  normal subgroup $M \cong \mathbb R^2$ of $G$ given in Theorem \ref{solvableonedimensional} (a) is the commutator subgroup
$G'= \{ \exp (t e_1 + u e_2);  t, u \in \mathbb R \}$
of $G$. Moreover, the inner mapping group $Inn(L)$ of $L$ is a $1$-dimensional non-normal subgroup of $G'$. Hence $Inn(L)$ must be the subgroup $H_4$ (see Theorem \ref{lefttranslation}). The normalizer of $H_4$ in $G$ is the group
$N=\{ \exp (t_1 e_1 + t_2 e_2+ t_3 e_3);  t_i \in \mathbb R \}$. As the direct product
$Z \times Inn(L)=G'$ we have a contradiction to Lemma \ref{niemenmaa}. \qed

\medskip
\noindent
Now we treat $4$-dimensional solvable Lie groups which are direct products.

\begin{Prop}
There exists no connected topological loop $L$ such that the multiplication group of $L$ is a $4$-dimensional solvable Lie group which is the direct product of proper connected Lie groups.
\end{Prop}
\Pro By Lemmata \ref{simplyconnecteduj}, \ref{fourdim} and \ref{simplyconnected} we may assume that the loop $L$ is homeomorphic to
$\mathbb R^3$. Every $4$-dimensional solvable decomposable Lie group has a $1$-dimensional normal subgroup $N$. As the orbit $N(e)$ is a $1$-dimensional  normal subloop of $L$ it follows from Lemma \ref{fourdim} that the group $Mult(L)$ is simply connected and its centre has dimension $\ge 1$. Hence $Mult(L)$ has the form $C \times S$, where $C$ is the group $\mathbb R$ and $S$ is a $3$-dimensional simply connected Lie group. The orbit $C(e)$ is a $1$-dimensional central subgroup of $L$ isomorphic to $\mathbb R$ (see Theorem 11 in \cite{albert1}). By Theorem \ref{solvableonedimensional} (a) there is a $2$-dimensional normal subgroup $M$ containing the group $C \cong \mathbb R$ and the commutator subgroup $Mult(L)'=S'$ of $Mult(L)$. Hence one has $\hbox{dim} \ Mult(L)'=1$. Then $Mult(L)$ is isomorphic either to $G_1=\mathbb R^2 \times {\mathcal L_2}$ or to $G_2=\mathbb R \times {\mathcal F}_3$, where ${\mathcal F_3}$ is the $3$-dimensional filiform Lie group.
Proposition 5.1 (i) in \cite{figula} shows that the group $G_2$ is not the multiplication group of a topological loop $L$ homeomorphic to
$\mathbb R^3$.
\newline
\noindent
Now we suppose that the group $Mult(L)$ is the group $G_1$ which is given on $\mathbb R^4$ by the multiplication
\begin{equation} \label{multuj}
g(x_1,x_2,x_3,x_4) g(y_1,y_2,y_3,y_4)=g(x_1+y_1, x_2+y_2, x_3+y_3, y_4+ x_4 e^{y_3}). \nonumber \end{equation}
Then the centre $Z$ of $G_1$ is the group $Z=\{ g(x,y,0,0), x,y \in \mathbb R \}$ and the commutator subgroup of $G_1$ is the group
$G_1'=\{ g(0,0,0,z), z \in \mathbb R \}$.
By Theorem 11 in \cite{albert1} the orbit $Z(e)$  is the centre of $L$ isomorphic to $\mathbb R^2$. Since the multiplication group $Mult(L/Z(e))$ of the factor loop $L/Z(e)$ is a factor group of $G_1$ (see Lemma \ref{brucklemma}) we get $L/Z(e)$ is the group $\mathbb R$
(see Theorem 18.18 in \cite{loops}). Hence there is a normal subgroup $P$ of $G_1$ such that $Z$ is a subgroup of $P$ and the factor group $G_1/P$ is isomorphic to the group  $Mult(L/Z(e)) \cong \mathbb R$ (Lemma \ref{brucklemma}). Then one has $G_1' < P$ and therefore $P=Z \times G_1'$. 
As $G_1/P$ acts sharply transitively on the orbits $Z(x)$, $x \in L$, the inner mapping group $Inn(L)$ of the loop $L$ is a $1$-dimensional subgroup of $P$ with $Co_{G_1}(Inn(L))=1$. The Lie algebra
${\bf g_1}$ of  $G_1$ has a basis $\{ e_1,e_2,e_3,e_4 \}$ with $[e_4,e_3]=e_4$. Hence the Lie algebra ${\bf p}$ of $P$ is given by ${\bf p}=\langle e_1, e_2, e_4 \rangle $ and we may choose $Inn(L)$ as the subgroup
$\exp t(e_4 + a e_1+ b e_2)$, $t \in \mathbb R$, with $a \neq 0$ or $b \neq 0$. Each automorphism
$\varphi $ of ${\bf g_1}$ has the form $\varphi (e_1)= k_1 e_1+ k_2 e_2$, $\varphi (e_2)= l_1 e_1+ l_2 e_2$,
$\varphi (e_4)=n e_4$, $\varphi (e_3)= a_1 e_1+ a_2 e_2+a_3 e_4 +e_3$ such that $(k_1 l_2- l_1 k_2)n \neq 0$, $k_i, l_i, n, a_j \in \mathbb R$, $i=1,2$, $j=1,2,3$. Then we can change $Inn(L)$ by an automorphism of $G_1$ such that $Inn(L)=\{ \exp t(e_4+ e_1), t \in \mathbb R \}=
\{g(u,0,0,u), u \in \mathbb R \}$.
\newline
\noindent
According to Lemma \ref{kepka} the group $G_1$ is isomorphic to the multiplication group $Mult(L)$ of a topological proper loop $L$ having the subgroup $Inn(L)$ as its inner mapping group precisely if there are two left transversals $A$ and $B$ to $Inn(L)$ in $G_1$ such that
$\{a^{-1} b^{-1} a b;\ a \in A, b \in B \}$ is contained in $Inn(L)$
and the set $\{ A, B \}$ generates the group $G_1$.
Arbitrary left transversals to the group $Inn(L)$ in $G_1$ are: $A=\{ g(x,y,z,f(x,y,z)); x,y,z \in \mathbb R \}$ and
$B=\{ g(k,l,m,h(k,l,m));\  k,l,m \in \mathbb R \}$,
where $f: \mathbb R^3 \to \mathbb R$, $h: \mathbb R^3 \to \mathbb R$ are continuous functions with $f(0,0,0)=h(0,0,0)=0$.
The products $a^{-1} b^{-1} a b$ with $a \in A$ and
$b \in B$ are elements of $Inn(L)$ if and only if the equation
$h(k,l,m) (1-e^z)= f(x,y,z) (1-e^m)$
holds for all $x,y,z,k,l,m \in \mathbb R$. Since the left hand side of the last equation does not depend on the variables
$x$ and $y$ and the right hand side is independent of $k$, $l$ we have $h(k,l,m)=h(m)$, $f(x,y,z)=f(z)$ and it follows that
$\frac{h(m)}{1-e^m}= \frac{f(z)}{1-e^z}=k$, where $k$ is a real constant. Then both sets $A$ and $B$ consist of the centre $Z$ of $G_1$ and the one-parameter subgroup $F= \{ g(0,0,z,k(1-e^z)), z \in \mathbb R \}$ with $Z \cap F=\{ 1 \}$. Hence $\{ A, B \}$ does not generate the group $G_1$. This  contradiction proves the assertion.   \qed

\begin{Prop} \label{4dimnoonedim} A $4$-dimensional connected Lie group having no  normal subgroup of dimension $1$ cannot be the multiplication group of a connected topological proper loop $L$.
\end{Prop}
\Pro  We may suppose that $L$ is homeomorphic to $\mathbb R^3$ (see Lemmata \ref{simplyconnecteduj}, \ref{fourdim} and \ref{simplyconnected}).
Any $4$-dimensional connected Lie group having no $1$-dimensional normal subgroup is locally isomorphic to the group $G$ given in Case 4.12 of \cite{ghanam}.
The Lie algebra ${\bf g}$ of $G$ is given by
$[e_1,e_3]=e_1,  \ [e_2,e_3]=e_2, \  [e_1,e_4]=-e_2, \ [e_2,e_4]=e_1$ (see ${\it g}_{4,10}$ in \cite{mubarakjzanov1}).

The commutator subgroup $G'$ of $G$ is the $2$-dimensional abelian normal subgroup $G'=\{ g(x,y,0,0), x,y \in \mathbb R \}$.
The orbit $G'(e)$ is a connected normal subloop of $L$ with dimension $1$ or $2$. As $G$ has discrete centre one has $\hbox{dim} G'(e)=2$ (see Lemma \ref{fourdim}). The multiplication group of the subloop $G'(e)$ is a subgroup of $G$ (see Lemma \ref{brucklemma}).   Then $G'(e)$ is isomorphic to $\mathbb R^2$ because none of the groups
$Mult({\mathcal L}_2)={\mathcal L}_2 \times {\mathcal L}_2$ and $Mult(L_{\mathcal F})={\mathcal F_n}$, $n \ge 4$, are contained in $G$. 
As the multiplication group of the factor loop $L/G'(e)$ is a factor group of $Mult(L)$ the loop $L/G'(e)$ is isomorphic to $\mathbb R$
(see Theorem 18.18 in \cite{loops}). Then there is a normal subgroup $K$ of $G$ such that $G/K$ is isomorphic to the multiplication group $Mult(L/G'(e)) \cong \mathbb R$ (cf. Lemma \ref{brucklemma}). Therefore the group $K$ has dimension $3$, it contains
the subgroup $G'$ and leaves every orbit $G'(x)$, $x \in L$, in $L$ invariant. Hence the Lie algebra ${\bf k}$ of $K$ has one of the following forms:
${\bf k_1}= \langle e_1, e_2, e_4+l e_3 \rangle$, $l \in \mathbb R$, ${\bf k_2}= \langle e_1, e_2, e_3 \rangle$.
The Lie group $K_1$ of ${\bf k_1}$ has no $1$-dimensional normal subgroup. For this reason $K_1$ cannot induce on the orbit $G'(e)$
a $2$-dimensional group.
Any $1$-dimensional normal subgroup $S$ of the Lie group $K_2$ of ${\bf k_2}$ is contained in the commutator subgroup $K_2'=G'$. Hence $K_2/S$ is isomorphic to ${\mathcal L}_2$. As $G'$ acts sharply transitively on $G'(e)$, for every element $s \in S \setminus \{ 1 \}$ one has $s(e) \neq e$ and $K_2$ cannot induce on the orbit $G'(e)$ a $2$-dimensional group.

Hence the group induced by $K_i$, $i=1,2$, on the orbit $G'(e)$ is isomorphic to $K_i$. Then $K_i$ induces a group
isomorphic to $K_i$ on every orbit $G'(x)$, $x \in L$. The same consideration as for the group $\Omega \cong {\mathcal L}_2$ discussed in the proof of Theorem \ref{solvableonedimensional} (a) is valid  for the groups $K_i$, $i=1,2$. 
Therefore the centre of $L$ would be at least $1$-dimensional and we have a contradiction to the fact that $G$ has discrete centre. \qed

\bigskip
\noindent
\centerline{\bf 5. Five-dimensional solvable indecomposable Lie groups}

\medskip
\noindent
There are 39 classes of $5$-dimensional solvable indecomposable Lie algebras (\cite{mubarakjzanov2}). Among them precisely the Lie algebras
${\it g}_{5,1}$ to ${\it g}_{5,6}$ are nilpotent.   The non-nilpotent Lie algebras have at most a $1$-dimensional centre.
In this section we prove that there does not exist $3$-dimensional connected topological loop $L$ such that the Lie algebra of the group $Mult(L)$ of $L$ is  a $5$-dimensional solvable non-nilpotent indecomposable Lie algebra.

\begin{Prop} \label{5dimdiscrete} There exists no $3$-dimensional connected topological proper loop $L$ such that the Lie algebra of its multiplication group is a $5$-dimensional solvable indecomposable Lie algebra with trivial centre.
\end{Prop}
\Pro We may assume that $L$ is homeomorphic to $\mathbb R^3$ (see Lemmata \ref{simplyconnecteduj} and \ref{simplyconnected}).
In \cite{mubarakjzanov2}  the $5$-dimensional solvable indecomposable Lie algebras ${\bf g}$ with trivial centre are  the Lie algebras
${\it g}_{5,7}$, ${\it g}_{5,9}$, the Lie algebras ${\it g}_{5,11}$ to ${\it g}_{5,13}$, the Lie algebras ${\it g}_{5,16}$ to ${\it g}_{5,18}$,
${\it g}_{5,21}$,
${\it g}_{5,23}$, ${\it g}_{5,24}$, ${\it g}_{5,27}$,  the Lie algebras ${\it g}_{5,31}$ to ${\it g}_{5,37}$,  the Lie algebras
${\it g}_{5,19}$, ${\it g}_{5,20}$ and ${\it g}_{5,28}$ in the case of that $\alpha \neq -1$,
${\it g}_{5,15}$ in the case of that $\gamma \neq 0$,
${\it g}_{5,25}$ in the case of that $\beta \neq 0$, $p \neq 0$, ${\it g}_{5,26}$ in the case of that $p \neq 0$, ${\it g}_{5,30}$ in the case of that
$h \neq -2$.
\newline
\noindent
All Lie algebras ${\bf g}$ from this list with exceptions of the Lie algebras ${\it g}_{5,17}$,  ${\it g}_{5,18}$ and ${\it g}_{5,33}$ have the $1$-dimensional ideal
${\bf n}_1=\langle e_1 \rangle $ such that the factor algebras ${\bf g}/ {\bf n}_1$ are not isomorphic to the Lie algebras of the groups
${\mathcal L}_2 \times {\mathcal L}_2$ or ${\mathcal F}_4$. As the centre of ${\bf g}$ is trivial  these Lie algebras cannot be the Lie algebras of the multiplication groups of $3$-dimensional topological loops (Theorem \ref{solvableonedimensional}).
The Lie algebra ${\it g}_{5,33}$ is defined by $[e_1,e_4]=e_1$, $[e_3,e_4]=\beta e_3$, $[e_2,e_5]=e_2$, $[e_3,e_5]=\gamma e_3$, where
$\gamma^2 + \beta^2 \neq 0$. The factor algebra ${\it g}_{5,33}/ \langle e_1 \rangle $, respectively ${\it g}_{5,33}/ \langle e_2 \rangle $ is isomorphic to the Lie algebra of ${\mathcal L}_2 \times {\mathcal L}_2$ precisely if $\gamma =0$, respectively $\beta =0$. But for $\gamma = \beta =0$ the Lie algebra ${\it g}_{5,33}$ is decomposable.
\newline
\noindent
Hence it remains to investigate the Lie algebras ${\it g}_{5,17}$ and ${\it g}_{5,18}$ which have no $1$-dimensional ideal. We denote by $G$ the Lie group of the Lie algebra ${\it g}_{5,17}$, respectively of ${\it g}_{5,18}$ and assume that $G$ is the multiplication group $Mult(L)$ of $L$. 
In both cases we consider the normal subgroup $N=\{ \exp (t_1 e_1+ t_2 e_2); t_i \in \mathbb R, i=1,2 \}$ of $G$.
\newline
\noindent 
First we suppose that the orbit $N(e)$ is a one-dimensional connected normal subloop of $L$. By Lemma \ref{brucklemma} the group $G$ has a connected normal subgroup $M$ containing the group $N$ such that the factor group $G/M$ is isomorphic to the multiplication group of the factor loop $L/N(e)$. Since $\hbox{dim} M \ge \hbox{dim} N=2$ the dimension of $G/M$ is $\le 3$. Hence by
Theorem \ref{solvableonedimensional} the factor group $G/M$ would be isomorphic to $\mathbb R^2$. As $G$ has discrete centre we have a contradiction to Theorem \ref{solvableonedimensional} (a).
\newline
\noindent
Therefore $N(e)$ is a two-dimensional connected normal subloop of $L$. The multiplication group $Mult(N(e))$
of $N(e)$ is a subgroup of $Mult(L)=G$. As none of the groups $Mult({\mathcal L_2})= {\mathcal L_2} \times {\mathcal L_2}$ and
$Mult(L_{\mathcal F})= \mathcal F_{n}$, $n \ge 4$, are subgroups of $G$ the normal subloop $N(e)$ is isomorphic to the group $\mathbb R^2$.
The multiplication group of the factor loop $L/N(e)$ is isomorphic to $\mathbb R$ (see Theorem 18.18 in \cite{loops}). There exists a normal subgroup $K$ of $G$ such that the factor group $G/K$ is isomorphic to $Mult(L/N(e)) \cong \mathbb R$ (see Lemma \ref{brucklemma}).
Hence $K$ contains the commutator subgroup $G'$ of $G$. Since $\hbox{dim} K= \hbox{dim} G'=4$ the group  $K$ coincides with the abelian group $G'$. Hence $K$ induces on the orbit $N(e)$ the group $\mathbb R^2$. The stabilizer $K_e$ of $e \in L$ in $K$ fixes every point on the orbit $N(e)=K(e)$. The inner mapping group $Inn(L)$ of $L$ is the group $K_e$. Hence $N(e)$ would be a $1$-dimensional central subgroup of $L$  which contradicts the fact that $G$ has discrete centre and the assertion follows.  \qed

\begin{Prop} \label{5dim1dimcentre} Let $L$ be a connected simply connected topological proper loop of dimension $3$ such that the Lie algebra of its multiplication group is a $5$-dimensional solvable non-nilpotent indecomposable Lie algebra having a $1$-dimensional centre. 
Then for the pair $({\bf g}, {\bf m})$ of the Lie algebras of the multiplication group $Mult(L)$ of $L$ and the abelian normal subgroup $M$ given in Theorem \ref{solvableonedimensional} (a) one of the following cases can occur:
\newline
\noindent
(a) The Lie algebra ${\bf g_1}$ is defined by $[e_2,e_3]=e_1$, $[e_2,e_5]=e_3$, $[e_4,e_5]=e_4$ and ${\bf m_1}= {\bf g_1}'$.
\newline
\noindent
(b) The Lie algebra ${\bf g_2}$ is defined by $[e_2,e_4]=e_1$, $[e_1,e_5]=e_1$, $[e_2,e_5]=e_2$, $[e_4,e_5]=e_3$ and
${\bf m_2}= {\bf g_2}'$.
\newline
\noindent
(c) The Lie algebra ${\bf g_3}$ is defined by $[e_1,e_4]=e_1$, $[e_2,e_5]=e_2$, $[e_4,e_5]=e_3$ and ${\bf m_3}= {\bf g_3}' $.
\newline
\noindent
(d) The Lie algebra ${\bf g_4}$ is defined by $[e_1,e_4]=e_1$, $[e_2,e_4]=e_2$, $[e_1,e_5]=-e_2$, $[e_2,e_5]=e_1$, $[e_4,e_5]=e_3$ and ${\bf m_4}= {\bf g_4}'  $.
\end{Prop}
\Pro By Lemma \ref{simplyconnected} the loop $L$ is homeomorphic to $\mathbb R^3$. According to  \cite{mubarakjzanov2} the $5$-dimensional solvable non-nilpotent indecomposable Lie algebras ${\bf g}$ with $1$-dimensional centre ${\zeta}$ are the Lie algebras ${\it g}_{5,8}$,
${\it g}_{5,10}$, ${\it g}_{5,14}$, ${\it g}_{5,22}$, ${\it g}_{5,29}$, ${\it g}_{5,38}$, ${\it g}_{5,39}$,
the Lie algebras ${\it g}_{5,19}$, ${\it g}_{5,20}$ and ${\it g}_{5,28}$ in the case of that $\alpha =-1$, ${\it g}_{5,15}$ in the case of that $\gamma =0$, ${\it g}_{5,25}$ in the case of that $\beta \neq 0$, $p=0$, ${\it g}_{5,26}$ in the case of that $p=0$, $\epsilon = \pm 1$ and
${\it g}_{5,30}$ in the case of that $h=-2$.
If ${\bf g}$ is the Lie algebra of the multiplication group $Mult(L)$ of $L$, then
the Lie group $Z= \exp {\zeta}$ is the centre of $Mult(L)$ and
 the orbit $Z(e)$, where $e$ is the identity element of $L$, is the $1$-dimensional centre of $L$ (see Theorem 11 in \cite{albert1}).
If $Mult(L)$ does not belong to the Lie algebra  ${\it g}_{5,38}$, then the factor algebras ${\bf g}/{ \zeta}$ are different from the Lie algebras of the Lie groups ${\mathcal L}_2 \times {\mathcal L}_2$ or ${\mathcal F}_4$.
Therefore the factor loop $L/Z(e)$ is isomorphic to $\mathbb R^2$ (cf. Theorem \ref{solvableonedimensional}). 
 The Lie algebra ${\it g}_{5,38}$ is defined by $[e_1,e_4]=e_1$, $[e_2,e_5]=e_2$, $[e_4,e_5]=e_3$.
As $S=\{ \exp (t e_1); t \in \mathbb R \}$ is a connected normal subgroup of the Lie group of ${\it g}_{5,38}$ the orbit $S(e)$ is a $1$-dimensional connected normal subloop of $L$. The factor algebra ${\it g}_{5,38}/ \langle e_1 \rangle $ is also different from the Lie algebras of the groups ${\mathcal L}_2 \times {\mathcal L}_2$ and ${\mathcal F}_4$ and the factor loop $L/S(e)$ is again isomorphic to $\mathbb R^2$. 
Hence the Lie algebra ${\bf g}$ of $Mult(L)$ has a $3$-dimensional abelian ideal ${\bf m}$ such that ${\bf m}$ contains the commutator ideal
${\bf g}'$ of ${\bf g}$ (cf. Theorem \ref{solvableonedimensional} (a)). The commutator ideal of the Lie algebras ${\it g}_{5,19}$, ${\it g}_{5,25}$, ${\it g}_{5,28}$ and ${\it g}_{5,30}$ has dimension $4$.
The commutator ideal of the Lie algebras ${\it g}_{5,20}$ and ${\it g}_{5,26}$ is non-abelian. Hence these Lie algebras cannot be the Lie algebras of the multiplication groups of $3$-dimensional topological loops.

For the Lie algebras ${\it g}_{5,8}$, ${\it g}_{5,10}$, ${\it g}_{5,14}$ and ${\it g}_{5,15}$ the commutator ideal ${\bf g}'$ of ${\bf g}$ is isomorphic to $\mathbb R^3$ and contains the centre of ${\bf g}$. Hence one has ${\bf m}= {\bf g}'$.  If
${\bf g}$ is the Lie algebra of the multiplication group of $L$, then
the Lie algebra ${\bf inn(L)}$ of the inner mapping group $Inn(L)$ of $L$ is a $2$-dimensional subalgebra of ${\bf m}$ containing no ideal $\neq 0$ of
${\bf g}$ (see Theorem \ref{solvableonedimensional} (a)). The direct sum of the centre
${\zeta }$ of ${\bf g}$ and the Lie algebra ${\bf inn(L)}$ coincides with ${\bf m}$. The Lie algebra ${\bf n}$ of the normalizer of $Inn(L)$ in the Lie group  of ${\bf g}$ is the $4$-dimensional abelian nilradical
${\it rad}= \langle e_1, e_2, e_3, e_4 \rangle$ of ${\bf g}$.  This contradiction to Lemma \ref{niemenmaa} yields that only the Lie algebras
${\it g}_{5,22}$, ${\it g}_{5,29}$, ${\it g}_{5,38}$ and ${\it g}_{5,39}$ can occur as the Lie algebras of the multiplication groups $Mult(L)$ of  $3$-dimensional topological loops $L$. The Lie algebra ${\it g}_{5,29}$ in \cite{mubarakjzanov2} is isomorphic to the Lie algebra given in assertion (b).
The ideal ${\bf m}$ of these Lie algebras is the commutator ideal and the assertion is proved. \qed

\medskip
\noindent
Now we exclude the Lie algebras in cases (a) to (d) of Proposition \ref{5dim1dimcentre}.

\begin{Prop} \label{elso}  There does not exist $3$-dimensional connected topological proper loop $L$ such that the Lie algebra ${\bf g}$ of the multiplication group of $L$ is one of the Lie algebras listed in cases (a) to (d) of Proposition \ref{5dim1dimcentre}.
\end{Prop}
\Pro  By Lemmata \ref{simplyconnecteduj} and  \ref{simplyconnected} we may assume that $L$ is homeomorphic to $\mathbb R^3$.
The linear representation of  the Lie group $G_i$ of ${\bf g_i}$ is:
For $i=1$
\begin{equation} g(x_1,y_1,z_1,q_1,w_1) g(x_2,y_2,z_2,q_2,w_2)= \nonumber \end{equation}
\begin{equation} g(x_1+w_1 y_2+\frac{w_1^2 z_2}{2}+x_2,y_1+w_1 z_2+y_2,z_1+ z_2, q_1+e^{z_1} q_2, w_1+w_2) \nonumber \end{equation}
for $i=2$
\begin{equation} g(q_1,x_1,y_1,z_1,w_1) g(q_2,x_2,y_2,z_2,w_2)= \nonumber \end{equation}
\begin{equation} g(q_1+e^{w_1} q_2+x_1 z_2,x_1+e^{w_1} x_2, y_1+w_1 z_2+y_2, z_1+z_2, w_1+w_2) \nonumber \end{equation}
for $i=3$
\begin{equation} g(q_1,x_1,y_1,z_1,w_1) g(q_2,x_2,y_2,z_2,w_2)= \nonumber \end{equation}
 \begin{equation} g(q_1+e^{z_1} q_2, x_1+ e^{w_1} x_2,y_1+ w_1 z_2 +y_2, z_1+z_2, w_1+ w_2). \nonumber \end{equation}
For $i=4$ the group $G_4$ is the linear group of matrices
\[ \left\{ g(x,y,q,w,z)= \left( \begin{array}{ccccc}
1 & x & y & -w & q \\
0 & e^w \cos z & e^w \sin z & 0 & 0 \\
0 & -e^w \sin z & e^w \cos z & 0 & 0 \\
0 & 0 & 0 & 1 & z \\
0 & 0 & 0 & 0 & 1 \end{array} \right), x,y,q,w,z \in \mathbb R \right\} \]
(cf. Cases 5.22, 5.29, 5.38, 5.39 in  \cite{ghanam}). First we determine which subgroups of the group $G_i$ can occur as the inner mapping group $Inn(L)_i$ of $L$. By Theorem \ref{solvableonedimensional} (a) the Lie algebra ${\bf inn(L)_i}$ of the inner mapping group $Inn(L)_i$ of $L$ 
is a $2$-dimensional subalgebra of the commutator ideal ${\bf m_i}={\bf g_i}'$ given in Proposition \ref{5dim1dimcentre} such that ${\bf inn(L)_i}$ does not contain any ideal $\neq \{ 0 \}$ of ${\bf g_i}$.

\noindent
As $\langle e_1 \rangle $ is the centre of ${\bf g_1}$ and $\langle e_4 \rangle $ is an ideal of ${\bf g_1}$ we may choose the Lie algebra
${\bf inn(L)_1}$  as follows
 ${\bf inn(L)_1}= \langle e_3+ a_1 e_1, e_4+ a_2 e_1 \rangle $, $a_1, a_2 \in \mathbb R$, $a_2 \neq 0$.
The automorphism group  of ${\bf g_1}$ consists of the following mappings
$\alpha (e_1)=c^2 e_1$, $\alpha (e_2)=b_1 e_1+ c e_2 + b_3 e_3$, $\alpha (e_3)=c f_3 e_1+ c e_3$, $\alpha (e_4)= d e_4$,
$\alpha (e_5)= f_1 e_1+ f_3 e_3+ f_4 e_4+ e_5$, where $c d \neq 0$, $b_1, b_3, f_1, f_3, f_4 \in \mathbb R$. Using an  automorphism of $G_1$ we may assume that
\begin{equation} \label{inn1} Inn(L)_1= \{ \exp (t e_3+ u(e_1+e_4)), t, u \in \mathbb R \}=\{ g(u,t,0,u,0), t, u \in \mathbb R \}. \nonumber \end{equation}
The centre of the Lie algebras ${\bf g_i}$, $i=2,3,4$, is $\langle e_3 \rangle $. Moreover, $\langle e_1 \rangle $, respectively $\langle e_1 \rangle $ and  $\langle e_2 \rangle $, respectively $\langle e_1, e_2 \rangle $ are ideals of ${\bf g_2}$, respectively ${\bf g_3}$, respectively ${\bf g_4}$. Hence 
we may choose ${\bf inn(L)_i}$, $i=2,3,4$, in the following way ${\bf inn(L)_i}=\langle e_1+ k_1 e_3, e_2+k_2 e_3 \rangle $, $k_1, k_2 \in \mathbb R$, such that for $i=2$ one has $k_1 \neq 0$, for $i=3$ we get $k_1 k_2 \neq 0$ and for $i=4$ at least one of the real parameters $k_1, k_2$ is different from $0$. 
For $k_1 k_2 \neq 0$ the automorphism $\alpha (e_1)=k_1 e_1$, $\alpha (e_2)=k_2 e_2$, $\alpha (e_3)=e_3$, $\alpha (e_4)=e_4$,
$\alpha (e_5)=e_5$ of ${\bf g_i}$, $i=2,3,4$, maps the Lie algebra ${\bf inn(L)_i}$ onto ${\bf inn(L)}_{2,1}= {\bf inn(L)_3}= {\bf inn(L)}_{4,1}= \langle e_1+ e_3, e_2+ e_3 \rangle $. For $k_2=0$ the automorphism $\gamma (e_1)= k_1 e_1$, $\gamma (e_2)=e_2$, $\gamma (e_3)=e_3$, $\gamma (e_4)=e_4$,
$\gamma (e_5)=e_5$ maps the subalgebra ${\bf inn(L)_i}$ onto ${\bf inn(L)}_{2,2}= {\bf inn(L)}_{4,3}= \langle e_1+ e_3, e_2 \rangle $. For $k_1=0$ the automorphism $\beta (e_1)= e_1$, $\beta (e_2)=k_2 e_2$, $\beta (e_3)=e_3$, $\beta (e_4)=e_4$,
$\beta (e_5)=e_5$ maps ${\bf inn(L)_i}$ onto ${\bf inn(L)}_{4,2}= \langle e_1, e_2+e_3 \rangle $. The corresponding Lie groups are 
$Inn(L)_{2,1}= Inn(L)_3= Inn(L)_{4,1}= \{ g(t_1,t_2,t_1+t_2,0,0), t_i \in \mathbb R, i=1,2 \}$, 
$Inn(L)_{2,2}= Inn(L)_{4,3}= \{ g(t_1,t_2,t_1,0,0), t_i \in \mathbb R, i=1,2 \}$, 
$Inn(L)_{4,2}= \{ g(t_1,t_2,t_2,0,0), t_i \in \mathbb R, i=1,2 \}$.

Arbitrary left transversals to the group $Inn(L)_i$ of $G_i$ are:
For $i=1$
\begin{equation} \label{transversal1} A_1=\{ g(k, f_1(k,l,m),l, f_2(k,l,m),m), k,l,m \in \mathbb R \}, \nonumber \end{equation}
\begin{equation} B_1=\{ g(u, g_1(u,v,w),v, g_2(u,v,w),w), u,v,w \in \mathbb R \}, \nonumber \end{equation}
for $i=2, 3, 4$
\begin{equation} \label{transversal2}
A=\{ g(k_1(k,l,m),k_2(k,l,m),k,l,m), k,l,m \in \mathbb R \} \nonumber \end{equation}
\begin{equation}
B=\{ g(h_1(u,v,w),h_2(u,v,w),u,v,w), u,v,w \in \mathbb R \}, \nonumber \end{equation}
where $f_i(k,l,m): \mathbb R^3 \to \mathbb R$, $k_i(k,l,m): \mathbb R^3 \to \mathbb R$, $g_i(u,v,w): \mathbb R^3 \to \mathbb R$, $h_i(u,v,w): \mathbb R^3 \to \mathbb R$, $i=1, 2$, are continuous functions with $f_i(0,0,0)=k_i(0,0,0)=g_i(0,0,0)=h_i(0,0,0)=0$.
We prove that none of the groups $G_i$, $i=1,2,3,4$, satisfies the condition that for all $a \in A_i$ and $b \in B_i$ one has 
$ a^{-1} b^{-1} a b \in Inn(L)_i$. It means that the groups $G_i$, $i=1,2,3,4$, are not multiplication groups of $L$ (cf. Lemma \ref{kepka}).

\noindent
The products $a^{-1} b^{-1} a b$ with $a=g(0,f_1(0,0,m),0,f_2(0,0,m),m) \in A_1$ and $b=g(0,g_1(0,v,0),v,g_2(0,v,0),0) \in B_1$
 are elements of $Inn(L)_1$ if and only if the equation
\begin{equation} \label{equequ3} f_2(0,0,m)= m \frac{g_1(0,v,0) e^v}{(1-e^{v})}- \frac{m^2 v e^v}{2(1-e^{v})}  \end{equation}
is satisfied for all $m,v \in \mathbb R$.
Since the left hand side of (\ref{equequ3}) depends only on the variable $m$ for all
$v \in \mathbb R \setminus \{ 0 \}$ the function $v \mapsto \frac{v e^v}{(1-e^{v})}$ must be constant which is a contradiction.

The products $a^{-1} b^{-1} a b$ with $a=g(k_1(0,0,m),k_2(0,0,m),0,0,m) \in A$, $b=g(h_1(0,v,0),h_2(0,v,0),0,v,0) \in B$ are contained in $Inn(L)_3$, respectively in $Inn(L)_{4,i}$, $i=1,2,3$,   
if and only if the equation
\begin{equation} \label{equel}   m=k_1(0,0,m) \frac{e^{-v}-1}{v}+ \frac{h_2(0,v,0)}{v}(1-e^{-m}),  \end{equation}
respectively for $i=1$ the equation 
\begin{equation} \label{equelelso}  -m= k_1(0,0,m) \frac{1- e^v}{v}+ k_2(0,0,m) \frac{1-e^v}{v}+ \frac{h_2(0,v,0)}{v}(\cos m - \sin m -1) + \nonumber \end{equation}
\begin{equation} \frac{h_1(0,v,0)}{v}(\cos m+ \sin m-1),  \end{equation}
respectively 
for $i=2$
\begin{equation} \label{equelmasodik} -m= k_2(0,0,m) \frac{1-e^v}{v} + \frac{h_1(0,v,0)}{v} \sin m + \frac{h_2(0,v,0)}{v}(\cos m-1),  \end{equation}
respectively 
for $i=3$
\begin{equation} \label{equelharmadik} -m= k_1(0,0,m) \frac{1-e^v}{v} + \frac{h_1(0,v,0)}{v}(\cos m-1)- \frac{h_2(0,v,0)}{v} \sin m   \end{equation}
holds for all $m,v \in \mathbb R$.
Since the left hand side of these equations depends only on the variable $m$ and the function $ v \mapsto \frac{1-e^{ \varepsilon v}}{v}$, where $\varepsilon =1$ or $-1$, is not constant we get  $k_j(0,0,m)=0$ and $h_j(0,v,0)=c_j v$, with $c_j \in \mathbb R$, $j=1,2$. Then equation (\ref{equel}), respectively (\ref{equelelso}), respectively (\ref{equelmasodik}), respectively (\ref{equelharmadik}) yields that for all $m \in \mathbb R$  the identity $m=c_2 (1-e^{-m})$, respectively 
$-m= c_1(\cos m+ \sin m-1)+ c_2(\cos m - \sin m -1)$,
respectively
$-m=c_1 \sin m + c_2(\cos m-1)$,
respectively
$-m=c_1(\cos m-1)- c_2 \sin m$
 is satisfied which is a contradiction.

The products
 $a^{-1} b^{-1} a b$ with $a=g(k_1(0,0,m),k_2(0,0,m),0,0,m) \in A$, $b=g(h_1(0,v,w),h_2(0,v,w),0,v,w) \in B$ are contained in $Inn(L)_{2,1}$, respectively in $Inn(L)_{2,2}$ if and only if the equation
\begin{equation} \label{haromequujequ} m v =  \end{equation}
\begin{equation} \frac{h_1(0,v,w)+h_2(0,v,w)}{e^w}(1- \frac{1}{e^m}) + \frac{k_1(0,0,m)}{e^m}(\frac{1}{e^w}-1)+\frac{k_2(0,0,m)}{e^m} (\frac{1+v}{e^w}-1), \nonumber \end{equation}
respectively
\begin{equation} \label{negyequujequ} m v = \frac{h_1(0,v,w)}{e^w}(1- \frac{1}{e^m}) + \frac{k_1(0,0,m)}{e^m}(\frac{1}{e^w}-1)+
\frac{v k_2(0,0,m)}{e^{m+w}} \end{equation}
is satisfied for all $m,v,w \in \mathbb R$. For $v=0$ equation (\ref{haromequujequ}), respectively (\ref{negyequujequ}) gives
$\frac{h_1(0,0,w)+h_2(0,0,w)}{1-e^w}= \frac{k_1(0,0,m)+k_2(0,0,m)}{1-e^m}= d$, respectively $\frac{h_1(0,0,w)}{1-e^w}= \frac{k_1(0,0,m)}{1-e^m}= d$ for a suitable constant $d \in \mathbb R$.
If $w=0$, then equation (\ref{haromequujequ}), respectively (\ref{negyequujequ}) yields
\begin{equation} \label{equujequujkilenc} v = \frac{h_1(0,v,0)+h_2(0,v,0)}{m e^m}(e^m-1) + \frac{k_2(0,0,m)}{m e^m} v, \end{equation}
respectively
\begin{equation} \label{equujequujtiz} v = \frac{h_1(0,v,0)}{m e^m}(e^m-1) + \frac{k_2(0,0,m)}{m e^m} v. \end{equation}
As the function $g: m \mapsto \frac{e^m-1}{e^m m}$ is not constant the right hand side of equation (\ref{equujequujkilenc}), respectively
(\ref{equujequujtiz}) is equal to $v$ precisely if $h_1(0,v,0)=-h_2(0,v,0)$ and $k_2(0,0,m)=m e^m$, respectively $h_1(0,v,0)=0$ and $k_2(0,0,m)=m e^m$.
Putting $k_1(0,0,m)=d(1-e^m)-k_2(0,0,m)$, $k_2(0,0,m)=m e^m$ into (\ref{haromequujequ}) and $k_1(0,0,m)=d(1-e^m)$, $k_2(0,0,m)=m e^m$ into  (\ref{negyequujequ}) we have
\begin{equation} \label{equujequtizenegy} v(e^w-1)= \frac{e^m-1}{m e^m}[h_1(0,v,w)+h_2(0,v,w)- d(1-e^w)], \end{equation}
respectively
\begin{equation} \label{equujequtizenketto} v(e^w-1)= \frac{e^m-1}{m e^m}[h_1(0,v,w)- d(1-e^w)]. \end{equation}
Since the left hand side of equations (\ref{equujequtizenegy}) and (\ref{equujequtizenketto}) depends only on the variables $v$ and $w$ and the function
$m \mapsto \frac{e^m-1}{m e^m}$ is not constant we get $h_1(0,v,w)+h_2(0,v,w)=d(1-e^w)$ in equation (\ref{equujequtizenegy}) and $h_1(0,v,w)=d(1-e^w)$ in equation (\ref{equujequtizenketto}). But then in both cases one has $v(e^w-1)=0$ for all $v,w \in \mathbb R$ which is a contradiction.
\qed

\bigskip
\noindent
\centerline{\bf 6. Three-dimensional topological loops having five-dimensional}
\centerline{\bf  solvable decomposable Lie groups as their multiplication groups}

\medskip
\noindent
We classify all $5$-dimensional connected solvable Lie groups which are direct products of proper connected subgroups and which are multiplication groups of $3$-dimensional connected simply connected topological loops $L$. Moreover, we determine the inner mapping groups of $L$.

\begin{Prop} \label{Propindecom} Let $L$ be a  connected simply connected  topological proper loop of dimension $3$ such that its multiplication group $Mult(L)$ is a $5$-dimensional solvable Lie group which is the direct product of connected subgroups. Then $L$ contains a central subgroup $C \cong \mathbb R$ such that the factor loop $L/C \cong \mathbb R^2$. Moreover:
\newline
\noindent
(I) If the centre of the group $Mult(L)$ has dimension $1$, then for the pair $({\bf mult(L)}, {\bf m})$ of the Lie algebras of $Mult(L)$ and the normal subgroup $M$ in Theorem \ref{solvableonedimensional} (a) one of the following cases occurs:
\newline
\noindent
(a) The group $Mult(L)_1$ is the group ${\mathcal F}_3 \times {\mathcal L}_2$. The Lie algebra ${\bf mult(L)_1}$ is defined by $[e_1,e_2]=e_3$, $[e_4,e_5]=e_4$  and
${\bf m_1}= \langle e_2, e_3, e_4 \rangle $.
\newline
\noindent
(b) The group $Mult(L)_2$ is the group ${\mathcal L}_2 \times {\mathcal L}_2 \times \mathbb R$. The Lie algebra ${\bf mult(L)_2}$ is defined by $[e_1,e_2]=e_1$, $[e_3,e_4]=e_3$, $[e_5,e_i]=0$ for all $i=1, \cdots ,4$, and
${\bf m_2}=\langle e_1, e_3, e_5 \rangle $.
\newline
\noindent
(c) The Lie algebra ${\bf mult(L)_3}$ is defined by $[e_2,e_3]=e_1$, $[e_1,e_4]=e_1$, $[e_2,e_4]=e_2$, $[e_5,e_i]=0$ for all $i=1, \cdots ,4$, and ${\bf m_3}=\langle e_1, e_2, e_5 \rangle $.
\newline
\noindent
(d) The Lie algebra ${\bf mult(L)_4}$ is defined by $[e_1,e_3]=e_1$, $[e_2,e_3]=e_2$, $[e_1,e_4]=-e_2$, $[e_2,e_4]=e_1$, $[e_5,e_i]=0$ for all $i=1, \cdots ,4$, and ${\bf m_4}=\langle e_1, e_2, e_5 \rangle $.
\newline
\noindent
(II) If $Mult(L)$ has $2$-dimensional centre, then it is either  the group ${\mathcal F}_4 \times \mathbb R$  or the direct product of the group $\mathbb R^2$ and a $3$-dimensional solvable Lie group $S$ having $2$-dimensional commutator subgroup. In the second case the Lie algebra ${\bf mult(L)}$ is the direct sum $\langle e_1, e_2, e_3 \rangle \oplus  \langle e_4, e_5 \rangle $, where
$\langle e_1, e_2, e_3 \rangle $ is the Lie algebra of $S$. The Lie algebra ${\bf m}$ has one of the following forms:
${\bf m}_{II,1}=\langle e_1, e_2, e_4 \rangle $, ${\bf m}_{II,2}=\langle e_1, e_2, e_5+k e_4 \rangle $, $k \in \mathbb R$.
\end{Prop}
\Pro The loop $L$ is homeomorphic to $\mathbb R^3$ (see Lemma \ref{simplyconnected}). We assume that the multiplication group $Mult(L)$ of $L$ is a $5$-dimensional decomposable solvable Lie group. Then for $Mult(L)$ we have the following possibilities: ${\mathcal L_2} \times \mathbb R^3$, ${\mathcal L_2} \times {\mathcal L_2} \times \mathbb R$, ${\mathcal L_2} \times S$, $\mathbb R^2 \times S$, $\mathbb R \times K$, where $S$ is a $3$-dimensional and $K$ is a $4$-dimensional solvable indecomposable Lie group.
All of these Lie groups have a normal subgroup $N \cong \mathbb R$ such that $Mult(L)/N$ is isomorphic neither to ${\mathcal L_2} \times {\mathcal L_2}$ nor to ${\mathcal F}_4$. Then the factor loop $L/N(e)$ is isomorphic to $\mathbb R^2$ (see Theorem \ref{solvableonedimensional}), the group $N(e)$ is central in $L$ and the first assertion is proved. Moreover, $Mult(L)$ is simply connected because it is a semidirect product of $\mathbb R^2$ with a normal subgroup $M \cong \mathbb R^3$ such that $M$ contains a $1$-dimensional central subgroup of $Mult(L)$ (cf. Theorem \ref{solvableonedimensional} (a)). 

Since $L$ is not associative, the centre $Z$ of $Mult(L)$ has dimension $1$ or $2$.
If $\hbox{dim} \ Z=1$, then $Mult(L)$ is either the group ${\mathcal F}_3 \times {\mathcal L}_2$ or
 the direct product
$K \times Z$, where $Z$ is the group $\mathbb R$ and $K$ is a $4$-dimensional solvable Lie group with discrete centre. 

If $Mult(L)= {\mathcal F}_3 \times {\mathcal L}_2$, then its Lie algebra ${\bf mult(L)}$ is given by $[e_1,e_2]=e_3$, $[e_4,e_5]=e_4$. The commutator ideal
${\bf mult(L)}'= \langle e_3, e_4 \rangle $ contains the centre
$\langle e_3 \rangle $ of ${\bf mult(L)}$. Since all $2$-dimensional subalgebras of the Lie algebra ${\bf  f_3}$ of ${\mathcal F}_3$ containing the centre of
${\bf f_3}$ can be mapped under an element of $Aut({\bf f_3})$ onto the subalgebra $\langle e_2, e_3 \rangle $
we may assume that the Lie algebra ${\bf m}$ of $M$ has the form as in case (a) of assertion (I).

If $Mult(L)= K \times Z$, then $Mult(L)$ has a normal subgroup $M \cong \mathbb R^3$ such that $M$ contains the commutator subgroup $Mult(L)'=K'$ and the centre $Z$ of $Mult(L)$.
Since there is no $4$-dimensional solvable Lie group with discrete centre and $1$-dimensional commutator subgroup, the dimension of $K'$ must be $2$. 
Hence the Lie algebra ${\bf k}$ of $K$ is one of the following: 
the Lie algebra of
${\mathcal L}_2 \times {\mathcal L}_2$ or  ${\it g}_{4,8}$ with $h=0$ or ${\it g}_{4,10}$ in \cite{mubarakjzanov1}, § 5. 
If ${\bf k}$ is the Lie algebra of ${\mathcal L}_2 \times {\mathcal L}_2$, then we get case (b) in assertion (I). If ${\bf k}$ is the Lie algebra
${\it g}_{4,8}$ with $h=0$, then we obtain case (c) of assertion (I). If ${\bf k}$ is the Lie algebra ${\it g}_{4,10}$, then we have case (d) in assertion (I).

Now we assume that $Mult(L)$ has a $2$-dimensional centre. If $Mult(L)$ is nilpotent, then it is the group ${\mathcal F}_4 \times \mathbb R$ and Proposition 5.1 of \cite{figula} proves the assertion. If $Mult(L)$ is not nilpotent, then it is either
the direct product
$K \times N$, where $N \cong \mathbb R$ and $K$ is a $4$-dimensional solvable non-nilpotent indecomposable Lie group with $1$-dimensional centre, or  the direct product $S \times R$, where $R \cong \mathbb R^2$ and $S$ is a $3$-dimensional solvable Lie group with discrete centre. 

If $Mult(L) = K \times N$, then the orbit $N(e)$ is a $1$-dimensional central subgroup of $L$ with $L/N(e) \cong \mathbb R^2$. 
Hence $Mult(L)$ has a normal subgroup $M \cong \mathbb R^3$ containing $N$ and the commutator subgroup $Mult(L)'=K'$ of $Mult(L)$. 
Among the $4$-dimensional solvable non-nilpotent Lie algebras only the Lie algebra ${\it g}_{4,3}$ has a $1$-dimensional centre and an abelian commutator subalgebra (cf. § 5 of \cite{mubarakjzanov1}).
If ${\bf k}$ is the Lie algebra ${\it g}_{4,3}$, then the Lie algebra ${\bf mult(L)}$ of $Mult(L)$ is defined by $[e_1,e_4]=e_1$, $[e_3,e_4]=e_2$, $[e_5,e_i]=0$ for all $i=1, \cdots ,4$, and the Lie algebra ${\bf m}$  of $M$ has the form $\langle e_1, e_2, e_5 \rangle $. The inner mapping group $Inn(L)$ of $L$ is a $2$-dimensional connected subgroup of $M$ such that $Co_{Mult(L)} (Inn(L))=1$. As $\langle e_2, e_5 \rangle $ is the centre of ${\bf mult(L)}$  the Lie algebra ${\bf inn(L)}$ of $Inn(L)$ has the form
${\bf inn(L)}= \langle e_2+ a_1 e_1, e_5+ a_2 e_1 \rangle$ with $a_1 a_2 \neq 0$. Then the Lie algebra $\langle e_1, e_2, e_3, e_5 \rangle$ of the normalizer $N_{Mult(L)} (Inn(L))$ is different from the Lie algebra $\langle e_1, e_2, e_5 \rangle$ of $Z \times Inn(L)$. This contradiction to Lemma \ref{niemenmaa} excludes the Lie algebra ${\it g}_{4,3}$.

If $Mult(L)=S \times R$, then the commutator ideal ${\bf i}=\langle e_1, e_2 \rangle $ of the Lie algebra ${\bf s}=\langle e_1, e_2, e_3 \rangle $ of $S$ is  commutative (see \cite{mubarakjzanov1}, § 4).
Let $N$ be a $1$-dimensional subgroup of the centre $R= \exp \{ a e_4+ b e_5, a, b \in \mathbb R \}$ of $Mult(L)$.
The Lie algebra ${\bf n}$ of $N$ has one of the following forms:
${\bf n}_1= \langle e_4 \rangle $, ${\bf n}_2= \langle e_5+k e_4 \rangle $, $k \in \mathbb R$.
As the Lie algebra ${\bf m}$ of the normal subgroup $M \cong \mathbb R^3$ is the direct sum
${\bf i} \oplus {\bf n}$, the form of
${\bf m}$ is given in assertion (II).
\qed

\begin{thm} \label{multinn} Let $L$ be a  connected simply connected topological proper loop of dimension $3$ such that its multiplication group is a $5$-dimensional solvable non-nilpotent Lie group which is the direct product of proper connected subgroups. Then the following Lie groups are the multiplication groups $Mult(L)$ and
the following subgroups are the inner mapping groups $Inn(L)$ of $L$:
\newline
\noindent
1) $Mult(L)_1$ is the Lie group ${\mathcal F}_3 \times {\mathcal L}_2$ the multiplication of which is given by
$g(x_1,x_2,x_3,x_4,x_5) g(y_1,y_2,y_3,y_4,y_5)=$
\begin{equation}  g(x_1+ y_1, x_2+y_2, x_3+ y_3- x_1 y_2, y_4+x_4 e^{y_5}, x_5+y_5). \nonumber \end{equation}
$Inn(L)_1$ is the following subgroup $\{ g(0,t,k,k,0); \ t,k \in \mathbb R \}$.
\newline
\noindent
2) $Mult(L)_2$ is the Lie group ${\mathcal L}_2 \times {\mathcal L}_2 \times \mathbb R$ which is represented on $\mathbb R^5$ by the multiplication  $g(x_1,x_2,x_3,x_4,x_5) g(y_1,y_2,y_3,y_4,y_5)= $
\begin{equation} g(y_1+ x_1 e^{y_2}, x_2+y_2, y_3+ x_3 e^{y_4}, x_4+y_4, x_5+y_5). \nonumber \end{equation}
$Inn(L)_2$ is the following subgroup $\{ g(t,0,k,0,t+k); \ t,k \in \mathbb R \}$.
\newline
\noindent
3) The multiplication of the group $Mult(L)_3$ is defined by
\begin{equation} \label{elem}
g(z_1,y_1,x_1,w_1,q_1) g(z_2,y_2,x_2,w_2,q_2)= \nonumber \end{equation}
\begin{equation} g(z_1+e^{w_1} z_2 -x_1 e^{w_1} y_2, y_1+ e^{w_1} y_2, x_1+x_2, w_1+w_2, q_1+q_2). \nonumber \end{equation}
$Inn(L)_3$ is one of the following groups: $Inn(L)_{3,1}=\{ g(z,y,0,0,z); z,y \in \mathbb R \}$,
$Inn(L)_{3,2}=\{ g(z,y,0,0,z+y); z,y \in \mathbb R \}$.
\newline
\noindent
4) The multiplication group $Mult(L)_4$ is the group  of matrices
\begin{equation} \label{elemek} \left\{ g(x,y,w,z,u)= \left( \begin{array}{cccc}
1 & x & y & u \\
0 & e^w \cos z & e^w \sin z & 0 \\
0 & -e^w \sin z &  e^w \cos z  & 0 \\
0 & 0 & 0 & 1 \end{array} \right), x,y,w,z,u \in \mathbb R \right\} \nonumber \end{equation}
(see Case 4.12 in \cite{ghanam}). Moreover, $Inn(L)_4$ is one of the following subgroups:
$Inn(L)_{4,1}= Inn(L)_{3,1}$,
$Inn(L)_{4,2}=\{ g(x,y,0,0,y); x,y \in \mathbb R \}$ 
$Inn(L)_{4,3}= Inn(L)_{3,2}$.
\newline
\noindent
5) The multiplication group $Mult(L)_5$ is the direct product of $\mathbb R^2$ and the connected Lie group of dimension $3$ having precisely one $1$-dimensional normal subgroup. The multiplication of $Mult(L)_5$ is given by
\begin{equation} g(x_1,x_2,x_3,x_4,x_5) g(y_1,y_2,y_3,y_4,y_5)= \nonumber \end{equation}
\begin{equation} g(y_1+x_1 e^{y_3}, y_2+x_2 e^{y_3}+x_1 y_3 e^{y_3},x_3+y_3,x_4+y_4,x_5+y_5). \nonumber \end{equation}
$Inn(L)_5$ is the following subgroup $\{ g(x,y,0,y,0); x,y \in \mathbb R \}$. 
\newline
\noindent
6) The elements of the multiplication group $Mult(L)_6$ can be written in the following form
\begin{equation} \label{3dimsincos} g(x,y,z,u,v)= \left( \begin{array}{ccccc}
1 & x & y & u & v \\
0 & e^{az} \cos z & e^{az} \sin z & 0 & 0 \\
0 & -e^{az} \sin z & e^{az} \cos z & 0 & 0 \\
0 & 0 & 0 & 1 & 0 \\
0 & 0 & 0 & 0 & 1 \end{array} \right), x,y,z,u,v \in \mathbb R, a > 0. \nonumber \end{equation}
$Inn(L)_6$ is one of the subgroups:
$Inn(L)_{6,1}=\{ g(x,y,0,x+y,0); x,y \in \mathbb R \}$, $Inn(L)_{6,2}=\{ g(x,y,0,x,0); x,y \in \mathbb R \}$,
$Inn(L)_{6,3}=Inn(L)_5$. 
\newline
\noindent
7) $Mult(L)_7$ is the direct product of $\mathbb R^2$ and the connected Lie group of dimension $3$ having precisely two $1$-dimensional normal subgroups. The group $Mult(L)_7$ is represented on $\mathbb R^5$ by the following multiplication
\begin{equation} \label{mult7} g(x_1,x_2,x_3,x_4,x_5) g(y_1,y_2,y_3,y_4,y_5)= \nonumber \end{equation}
\begin{equation} g(y_1+x_1 e^{a y_3}, y_2+x_2 e^{b y_3},x_3+y_3,x_4+y_4,x_5+y_5),   \end{equation}
with fixed but different numbers $a, b \in \mathbb R \backslash \{ 0 \}$.
\newline
\noindent
8) $Mult(L)_8$ is the direct product of $\mathbb R^2$ and the  connected Lie group of dimension $3$ having infinitely many $1$-dimensional normal subgroups. The multiplication of $Mult(L)_8$  is given by (\ref{mult7})
with $a=b \in \mathbb R \backslash \{ 0 \}$.
\newline
\noindent
The inner mapping group $Inn(L)_i$, $i=7,8$, is the group $Inn(L)_{6,1}$.
\end{thm}
\Pro By Lemma \ref{simplyconnected} the loop $L$ is homeomorphic to $\mathbb R^3$.
For $i=1,2,3,4$, the Lie algebras ${\bf mult(L)}_i$ of the groups $Mult(L)_i$ and the ideals ${\bf m}_i$ of ${\bf mult(L)}_i$ are given in Proposition
\ref{Propindecom}, (I) cases (a) to (d).  The Lie algebra ${\bf inn(L)}_i$ of the inner mapping group $Inn(L)_i$ of $L$ is a $2$-dimensional subalgebra of ${\bf m}_i$  containing no ideal $\neq \{0 \}$ of ${\bf mult(L)}_i$, $i=1,2,3,4$. For $i=1$ the Lie algebra ${\bf inn(L)}$ has the form ${\bf inn(L)}_{b_1,b_2}= \langle e_2+ b_1 e_3, e_4+ b_2 e_3 \rangle $, $b_1,b_2 \in \mathbb R$, $b_2 \neq 0$. The automorphism $\beta (e_1)=e_1$,
$\beta (e_2)=e_2-b_1 e_3$, $\beta (e_3)=e_3$, $\beta (e_4)= b_2 e_4$, $\beta (e_5)=e_5$
maps ${\bf inn(L)}_{b_1,b_2}$ onto ${\bf inn(L)}_1=\langle e_2, e_4+e_3 \rangle $.
The corresponding group $Inn(L)_1$ is given in assertion 1).
\newline
\noindent
As $\langle e_5 \rangle $ is the centre of ${\bf mult(L)}_2$ the Lie algebra ${\bf inn(L)_2}$ has the form ${\bf inn(L)}_{a_1,a_2}=\langle e_1+ a_1 e_5, e_3+ a_2 e_5 \rangle $, $a_1, a_2 \in \mathbb R$ with
$a_1 a_2 \neq 0$. Using the automorphism $\alpha (e_1)=a_1 e_1$, $\alpha (e_2)=e_2$, $\alpha (e_3)= a_2 e_3$, $\alpha (e_4)=e_4$, $\alpha (e_5)=e_5$ of 
${\bf mult(L)}_2$ the Lie algebra ${\bf inn(L)}_{a_1,a_2}$ is reduced to ${\bf inn(L)}_2=\langle e_1+e_5, e_3+e_5 \rangle $. The corresponding group $Inn(L)_2$ is given in assertion 2).
\newline
\noindent
As $\langle e_5 \rangle $ is the centre  and $\langle e_1, e_2 \rangle $ is the commutator ideal of
${\bf mult(L)_i}$ for $i=3,4$, we can write ${\bf inn(L)_i}$ in the form
${\bf inn(L)}_{k_1,k_2}= \langle e_1+ k_1 e_5, e_2+ k_2 e_5 \rangle $, $k_1$, $k_2 \in \mathbb R$. For $i=3$ one has $k_1 \neq 0$ and for $i=4$ at least one of the parameters $k_1, k_2$ is different from $0$. 
Similarly to the automorphism $\alpha $ of ${\bf mult(L)}_2$ we can find suitable automorphisms of ${\bf mult(L)_i}$, $i=3,4$, which map the Lie algebra ${\bf inn(L)}_{k_1,0}$ onto ${\bf inn(L)}_{3,1}={\bf inn(L)}_{4,1}= \langle e_1+e_5, e_2 \rangle $, the Lie algebra 
${\bf inn(L)}_{0,k_2}$ onto ${\bf inn(L)}_{4,2}= \langle e_1, e_2+ e_5 \rangle $ and the Lie algebra ${\bf inn(L)}_{k_1,k_2}$, $k_1 k_2 \neq 0$, onto ${\bf inn(L)}_{3,2}={\bf inn(L)}_{4,3}=\langle e_1+ e_5, e_2+ e_5 \rangle $. The corresponding Lie groups are the groups $Inn(L)_{3,1}=Inn(L)_{4,1}$, $Inn(L)_{3,2}=Inn(L)_{4,3}$, $Inn(L)_{4,2}$ given in assertions 3) and 4).

\noindent
The sets
$A_1=\{ g(x,e^z-1,y,0,z); x,y,z \in \mathbb R \}$ and $B_1=\{ g(n,0,l,-n,m); l, \\ m,n \in \mathbb R \}$ are $Inn(L)_1$-connected left transversals in $Mult(L)_1$.  The sets $A_2=\{ g(2-e^{x_2}-e^{x_4}, x_2, 0, x_4, x_5+2-e^{x_2}-e^{x_4}); x_2,x_4,x_5 \in \mathbb R \}$ and $B_2=\{ g(1-e^{y_2}, y_2, 1-e^{y_2}, y_4, y_5); y_2,y_4,y_5 \in \mathbb R \}$
are $Inn(L)_2$-connected transversals in $Mult(L)_2$. The sets $A_3=\{ g((e^w-1)(x+2)-x,1-e^w,x,w,q); x,w,q \in \mathbb R \}$ and $B_3=\{ g((2-e^l)k,e^l-1,k,l,m); k,l,m \in \mathbb R \}$, respectively
the sets $B_3$ and $C_3=\{ g(x(e^w-2),1-e^w,x,w,q); x,w,q \in \mathbb R \}$
 are $Inn(L)_{3,1}$-, respectively $Inn(L)_{3,2}$-connected transversals in $Mult(L)_3$. The set $A_4=B_4=\{g(1- e^u \cos v, - e^u \sin v, u,v,w); u,v,w \in \mathbb R \}$ is a left transversal to the subgroups $Inn(L)_{4,i}$ for every $i=1,2,3$ in $Mult(L)_4$. Moreover, the sets $\{A_i, B_i \}$ for all $i=1,2,3,4$ as well as $\{B_3, C_3 \}$ generate the group $Mult(L)_i$.  This proves assertions 1) to 4) (cf. Lemma \ref{kepka}).

\noindent
The Lie algebra ${\bf mult(L)_5}$ of the group $Mult(L)_5$ in assertion 5) is defined by $[e_1,e_3]=p e_1-e_2$, $[e_2,e_3]=e_1+ p e_2$, $[e_4,e_i]=[e_5,e_i]=[e_4,e_5]=0$, $i=1,2,3$, $p > 0$ (see ${\it g}_{3,5}$ in \cite{mubarakjzanov1}, § 4).
The Lie algebra ${\bf mult(L)_6}$ of the group $Mult(L)_6$ in assertion 6) is given by $[e_2,e_3]=e_2$, $[e_1,e_3]=e_1+e_2$,
$[e_1,e_2]=[e_4,e_5]=[e_4,e_i]= [e_5,e_i]=0$,  $i=1,2,3$ (see \cite{loops}, Lemma 23.16).
The Lie algebra ${\bf mult(L)_7}$ of the group $Mult(L)_7$ in assertion 7) is  defined by  $[e_1,e_3]=a e_1$, $[e_2,e_3]=b e_2$, $[e_1,e_2]=[e_4,e_i]= [e_5,e_i]=[e_4,e_5]=0$,  $i=1,2,3$, where $a \neq b \in \mathbb R \backslash \{ 0 \}$.  For $a=b$ we get the Lie algebra
${\bf mult(L)_8}$ of the group $Mult(L)_8$ in assertion 8) (see \cite{loops}, Section 23.1).
\newline
\noindent 
For $i=5,6,7,8$, the Lie algebra ${\bf inn(L)_i}$ of the inner mapping group $Inn(L)_i$  of $L$
is a $2$-dimensional subalgebra of ${\bf m}_{II,j}$, $j=1,2$,  given in Proposition \ref{Propindecom} (II) containing no ideal $\neq 0$ of ${\bf mult(L)_i}$.
The Lie algebra ${\bf inn(L)_i}$ has one of the following forms:
${\bf inn(L)}_{a_1,a_2}= \langle e_1+ a_1 e_4, e_2+ a_2 e_4 \rangle $, $a_1, a_2 \in \mathbb R$ and ${\bf inn(L)}_{b_1,b_2}= \langle e_1+ b_1(e_5+ k e_4), e_2+ b_2(e_5+ k e_4) \rangle $,  $b_1, b_2, k \in \mathbb R$, such that for $i=5$ one has $a_2 b_2 \neq 0$, for $i=6$ at least one of the parameters $a_1, a_2$, respectively $b_1, b_2$ is different from $0$, for $i=7, 8$, one has $ a_1 a_2 b_1 b_2 \neq 0$.
\newline
\noindent
For $i=5$ using the automorphism $\alpha (e_1)=e_1+ \frac{a_1}{a_2} e_2$, $\alpha (e_2)=e_2$, $\alpha (e_3)=e_3$, 
$\alpha (e_4)=\frac{1}{a_2} e_4$, 
$\alpha (e_5)=e_5$, respectively  $\beta (e_1)=e_1+ \frac{b_1}{b_2} e_2$, $\beta (e_2)=e_2$, $\beta (e_3)=e_3$, $\beta (e_4)=e_4+e_5$,
$\beta (e_5)=\left( \frac{1}{b_2} -k \right) e_4- k e_5$ of ${\bf mult(L)_5}$ we can change ${\bf inn(L)}_{a_1,a_2}$, respectively ${\bf inn(L)}_{b_1,b_2}$  onto the Lie algebra ${\bf inn(L)_5}= \langle e_1, e_2+ e_4 \rangle $. 
\newline
\noindent 
For $i=6,7,8$ the automorphism $\gamma (e_1)=a_1 e_1$, $\gamma (e_2)=a_2 e_2$, $\gamma (e_3)=e_3$, $\gamma (e_4)=e_4$,
$\gamma (e_5)=e_5$, respectively $\delta (e_1)=b_1 e_1$, $\delta (e_2)=b_2 e_2$, $\delta (e_3)=e_3$, $\delta (e_4)=e_4+ e_5$,
$\delta (e_5)=(1-k)e_4 -k e_5$ of ${\bf mult(L)_i}$ maps the Lie algebra ${\bf inn(L)}_{a_1,a_2}$, respectively ${\bf inn(L)}_{b_1,b_2}$  onto 
${\bf inn(L)}_{6,1}={\bf inn(L)}_{7}={\bf inn(L)}_{8}= \langle e_1+e_4, e_2+e_4 \rangle $.
The automorphism $\gamma $, respectively $\delta $ of ${\bf mult(L)}_6$ with $a_2=1=b_2$ maps the Lie algebra ${\bf inn(L)}_{a_1,0}$, respectively 
${\bf inn(L)}_{b_1,0}$  onto ${\bf inn(L)}_{6,2}= \langle e_1+ e_4, e_2 \rangle $. The automorphism $\gamma $, respectively $\delta $ of 
${\bf mult(L)}_6$ with $a_1=1=b_1$ maps ${\bf inn(L)}_{0,a_2}$, respectively
${\bf inn(L)}_{0,b_2}$  onto ${\bf inn(L)}_{6,3}= {\bf inn(L)}_5$. The corresponding Lie groups $Inn(L)_5=Inn(L)_{6,3}$, $Inn(L)_{6,2}$, $Inn(L)_{6,1}=Inn(L)_7= Inn(L)_8$ are given in assertions 5) to 8). 
\newline
\noindent 
The sets $A_5=\{ g(0,1-e^{k_1}(1+k_1),k_1,k_2+1-e^{k_1}(1+k_1),k_3); k_i \in \mathbb R, i=1,2,3 \}$ and
$B_5=\{ g(1-e^{l_1},1-e^{l_1},l_1,l_2,l_3); l_i \in \mathbb R, i=1,2,3 \}$ are $Inn(L)_5$-connected transversals in $Mult(L)_5$. 
\newline
\noindent 
The set 
$A_6=B_6=\{ g(1+e^{a k_1}(\sin k_1- \cos k_1), 1-e^{a k_1}(\sin k_1+ \cos k_1),k_1,k_2, \\ k_3); k_i \in \mathbb R \}$ is for every $i=1,2,3$, a left transversal to $Inn(L)_{6,i}$ in $Mult(L)_6$. The set $A_7=B_7= \{ g(2-e^{b k_1}-e^{a k_1},2-e^{b k_1}-e^{a k_1},k_1,k_2,k_3); k_i \in \mathbb R, i=1,2,3 \}$ is a left transversal to $Inn(L)_7$ in $Mult(L)_7$. The set $A_8=B_8=\{ g(1-e^{a k_1}-k_1,k_1, k_1,k_2,k_3); k_i \in \mathbb R, i=1,2,3 \}$ is  a left transversal to $Inn(L)_8$ in $Mult(L)_8$. 
Since Lemma \ref{kepka} is satisfied for all these transversals,  assertions 5) to 8) is proved.  
\qed

\medskip
\noindent
By the previous theorem  only a classification of connected simply connected $5$-dimensional solvable Lie groups which are the multiplication groups of  connected topological loops $L$ with dimension $3$ is given. The next proposition shows that Lie groups which cannot be the multiplication groups of $L$  can have universal coverings which are multiplication groups of $L$.

\begin{Prop} \label{sincos3dim} The direct product $G$ of $\mathbb R^2$ and the connected component of the euclidean motion group  of $\mathbb R^2$ cannot be the multiplication group of a $3$-dimensional topological loop $L$.
\end{Prop}
\Pro The group $G$ is represented in case 6) of Theorem  \ref{multinn} such that $a=0$.
The subgroups of $G$ which can occur as the inner mapping group of $L$ are also listed in case 6) of Theorem  \ref{multinn}.
Arbitrary left transversals to $Inn(L)_{6,i}$, $i=1,2,3$, are  $A=\{ g(f_1(k_1,k_2,k_3), f_2(k_1,k_2,k_3),k_1,k_2,k_3); k_i \in \mathbb R \}$ and
$B=\{ g(h_1(l_1,l_2,l_3), h_2(l_1,l_2,l_3),l_1,l_2,l_3); l_i \in \mathbb R \}$
 such that for the continuous functions
$f_j(k_1,k_2,k_3): \mathbb R^3 \to \mathbb R$,
$h_j(l_1,l_2,l_3): \mathbb R^3 \to \mathbb R$, $j=1,2$, one has $f_j(0,0,0)=h_j(0,0,0)=0$. The set $\{a^{-1} b^{-1} a b, \ a \in A, b \in B \}$ is contained in  $Inn(L)_{6,i}$ if and only if for $i=1$
\begin{equation} \label{equsin3} h_1(l_1,l_2,l_3)(1-\cos{k_1}- \sin{k_1})+ h_2(l_1,l_2,l_3)(1+ \sin{k_1}-\cos{k_1})= \nonumber \end{equation}
\begin{equation} f_1(k_1,k_2,k_3)(1- \cos{l_1}- \sin{l_1})+ f_2(k_1,k_2,k_3) (1+ \sin{l_1}-\cos{l_1}) \end{equation}
for $i=2$
\begin{equation} \label{equsin1} h_1(l_1,l_2,l_3)(1- \cos{k_1})+ h_2(l_1,l_2,l_3) \sin{k_1}= \nonumber \end{equation}
\begin{equation} f_1(k_1,k_2,k_3)(1- \cos{l_1})+ f_2(k_1,k_2,k_3) \sin{l_1} \end{equation}
for $i=3$
\begin{equation} \label{equsin2} h_2(l_1,l_2,l_3)(1- \cos{k_1})- h_1(l_1,l_2,l_3) \sin{k_1}= \nonumber \end{equation}
\begin{equation} f_2(k_1,k_2,k_3)(1- \cos{l_1})- f_1(k_1,k_2,k_3) \sin{l_1} \end{equation}
holds for all $k_1,k_2,k_3,l_1,l_2,l_3 \in \mathbb R$. As the right hand side of equations (\ref{equsin3}), (\ref{equsin1}) and  (\ref{equsin2})  does not depend on the variables $l_2$, $l_3$ and the left hand side of (\ref{equsin3}), (\ref{equsin1}) and (\ref{equsin2}) is independent of $k_2$, $k_3$ we get $h_j(l_1,l_2,l_3)=h_j(l_1)$ and $f_j(k_1,k_2,k_3)=f_j(k_1)$ for all $j=1,2$. In this case the function $h_j(l_1)$, respectively $f_j(k_1)$, $j=1,2$, has the form $a_{1,j}(1- \cos l_1)+a_{2,j} \sin l_1$, respectively $b_{1,j}(1-\cos k_1)+b_{2,j} \sin k_1$, where $a_{1,j}, a_{2,j}, b_{1,j}, b_{2,j} \in \mathbb R$. Then the set $A \cup B$ does not generate $G$.
This contradiction to Lemma \ref{kepka} yields the assertion. \qed

\medskip
\noindent
Acknowledgment

This paper was supported by the Hungarian Scientific Research Fund (OTKA) Grant PD 77392, by the EEA and Norway Grants (Zolt\'an Magyary Higher Education Public Foundation) and by the J\'anos Bolyai Research Fellowship.


\begin{thebibliography}{37}

\bibitem{albert1} Albert, A. A. (1943). Quasigroups I. Trans. Amer. Math. Soc. 54:507-519.


\bibitem{bruck} Bruck, R. H. (1958). A Survey of Binary Systems. Berlin-G\"ottingen-Heidelberg: Springer-Verlag.

\bibitem{niem3} Cs\"org\"o,  P., Niemenmaa,  M. (2002). On connected transversals to nonabelian subgroups. European J. Combin. 23:179-185.

\bibitem{figula0} Figula, \'A. (2009). The multiplication groups of 2-dimensional topological loops. J. Group Theory  12:419-429.

\bibitem{figula} Figula, \'A. (2011). On the multiplication groups of 3-dimensional topological loops. J. Lie Theory 21:385-415.

\bibitem{figula2} Figula, \'A. Three-dimensional loops as sections in a four-dimensional solvable Lie group.  Proceedings of Iscia Group Theory, Ischia, Italy, Apr. 14-17, 2010; Bianchi, M., Longobardi, P., Maj, M., Eds.;  World Scientific Publishing Company, 2011.


\bibitem{ghanam} Ghanam, R., Strugar, I.,  Thompson, G. (2005). Matrix Representation for Low Dimensional Lie Algebras. Extracta Math. 20:151-184.

\bibitem{gorbatsevich} Gorbatsevich, V. V. (1977).  On three-dimensional homogeneous spaces.  Sib. Math. J. 18:200-210.

\bibitem{hofmann} Hofmann, K. H., Strambach, K. (1990). Topological and analytic loops. In: Chein, O., et al., ed. Quasigroups and Loops: Theory and Applications. Berlin: Heldermann-Verlag, pp. 205-262.



\bibitem{mostow} Mostow, G. D. (1954). Factor spaces of solvable groups. Ann. of Math. 60:1-27.


\bibitem{mubarakjzanov1} Mubarakzjanov, G. M. (1963). On solvable Lie algebras. (Russian) Izv. Vyssh. Uchebn. Zaved. Matematika 1:114-123.

\bibitem{mubarakjzanov2} Mubarakzjanov, G. M. (1963). Classification of real structures of Lie algebras of fifth order. (Russian) Izv. Vyssh. Uchebn. Zaved. Matematika 3:99-106.


\bibitem{niem4} Niemenmaa, M. (1995). On loops which have dihedral 2-groups as inner mapping groups. Bull. Austral. Math. Soc. 52:153-160.

\bibitem{kepka} Niemenmaa, M.,  Kepka, T. (1990). On Multiplication Groups of Loops. J.  Algebra 135:112-122.

\bibitem{niem1} Niemenmaa, M.,  Kepka, T. (1994). On connected transversals to abelian subgroups. Bull. Austral. Math. Soc.  49:121-128.

\bibitem{niem2} Niemenmaa, M.,  Vesanen, A. (1994). On connected transversals in the projective special linear group $PSL(2,7)$. J. Algebra 166:  455-460.

\bibitem{loops} Nagy,  P. T., Strambach, K. (2002). Loops in Group Theory
and Lie Theory.  Berlin-New York: Walter de Gruyter.


\bibitem{salzmann}  Salzmann, H. (1963). Zur Klassifikation topologischer Ebenen. Math. Ann. 150:226-241.

\bibitem{vesanen} Vesanen, A. (1996). Solvable loops and groups. J. Algebra 180:862-876.


\end{thebibliography}
\end{document}